\documentclass[11pt]{article}

\usepackage[utf8]{inputenc}
\usepackage[T1]{fontenc}
\usepackage{amssymb}
\usepackage{amsmath}
\usepackage{siunitx}
\usepackage{graphicx}
\usepackage{natbib}
\usepackage{caption}
\usepackage{subcaption}
\usepackage{authblk}
\usepackage{hyperref}

\hypersetup{breaklinks=true}

\graphicspath{{FiguresFullSWOF/}}

\title{\vspace{-2\baselineskip}FullSWOF: A free software package for the simulation of shallow water flows}

\author[1]{Olivier Delestre\thanks{Corresponding author. Olivier.Delestre@univ-cotedazur.fr}}
\author[2]{Fr\'ed\'eric Darboux}
\author[3]{Fran\c{c}ois James}
\author[3]{Carine Lucas}
\author[3]{Christian Laguerre}
\author[3]{St\'ephane Cordier\small}
\affil[1]{ Lab. J.A. Dieudonn\'e \& EPU Nice Sophia, Univ. Nice, France.}
\affil[2]{Inra, UR0272, UR Science du sol, Centre de recherche Val de Loire, CS 40001, F-45075 Orl\'eans Cedex 2, France.}
\affil[3]{MAPMO, UMR CNRS 7349, F\'ed\'eration Denis Poisson, FR CNRS 2964, Universit\'e d'Orl\'eans, F-45067 Orl\'eans cedex 02, France.}
\normalsize

\date{}

\begin{document}

\maketitle
\vspace{-0.5cm}

\section*{Abstract}
Numerical simulations of flows are required in numerous applications and
are typically performed using shallow-water equations. 
Here, we describe the FullSWOF software, which is based on up-to-date
finite volume methods and on well-balanced schemes
used to solve these types of equations. The software consists of a set
of open-source C++ codes and is freely available to the community,
easy to use, and open for further development. 
Several features make FullSWOF particularly suitable 
for applications in hydrology: small water depths and wet-dry transitions are robustly addressed, rainfall and infiltration
are incorporated, and data from grid-based digital topographies can be used directly.
A detailed mathematical description is given here, and
the capabilities of FullSWOF are illustrated based on analytic
solutions and on datasets compiled from real cases.
The codes, available in 1D and 2D versions,
 have been validated on a large set of benchmark cases, which are available together with the download information 
and documentation at \url{https://www.idpoisson.fr/fullswof/}.

\paragraph{Keywords :}
Overland flow; modeling; shallow-water equations; software; open source

\section{Introduction}\label{Sec:Intro}

The simulation of geophysical surface flows is required when analyzing
a large variety of natural and man-made situations. 
For such purposes, several mathematical models are available. 
More complex models utilize the full system of Navier-Stokes equations; these equations can be solved by 
Gerris, for example \citep{popinet11}. However, the parametrization of
the model is very intricate, and people typically prefer simplified
models that are chosen based on the characteristics of the simulation. 
For avalanche modeling (on large slopes), a simplification of the Navier-Stokes equations leads to the 
Savage-Hutter system \citep{Savage91}. In the case of shallow flows,
if the vertical motion of the fluid can be neglected, 
the shallow-water equations can be used. Additional simplifications may be
needed. For example, if the shallow-water equations continue to be too complex, they can be replaced by the diffusive wave or kinematic wave equations \citep{Moussa00, Novak10}.
In this hierarchy of mathematical models, the present paper considers the shallow-water equations and presents the FullSWOF program. 
The shallow-water equations were first formulated by Saint-Venant \citep{saintvenant71} 
in the study of floods and tides. These equations are now classical
equations in the fields of flood forecasting, pollutant transport, dam breaks, tsunamis, soil erosion by overland flow, etc.
Because no explicit solutions to the shallow-water equations are
known for general situations,
efficient and robust numerical simulations are required. 

Over the last forty years, numerous codes that make use of various methods have been developed in this field. The MacCormack scheme has been widely used for scientific purposes \citep[e.g.,][]{Zhang89, Esteves00, Fiedler00}.
Although it is relatively easy to program and performs computations quickly,
this scheme neither guarantees the positivity of water depths at the wet-dry transitions
nor preserves steady states, i.e., it is not well balanced
\citep{Lee10}, which requires some work if these issues are
to be addressed
\citep[e.g.,][]{Esteves00, Fiedler00}.
In the various industrial codes used in engineering (e.g., CANOE
\citep{Tanguy13}, HEC-RAS \citep{HECRAS-Hyd11}, ISIS \citep{ISIS}, and
MIKE11 \citep{MIKE11}), the flow equations, namely, the shallow-water equations,
are often solved in a non-conservative form \citep{Novak10} with
either the Preissmann scheme
or the Abbott-Ionescu scheme, leading to inaccurate calculations for transcritical flows
and hydraulic jumps. 

While all these pieces of software have been used for scientific research,
most of their source codes have not been made available to the community, raising a major concern about
research reproducibility. Reproducibility is a key component of the scientific method and 
has received increased interest in recent years in the computational modeling community.
Numerical simulations rarely ensure this essential property, leading to a low 
confidence in scientific results and undermining the advancement of knowledge \citep{Claerbout92, Stodden13}. While the free availability of a source code is not a
sufficient condition to make numerical results reproducible, it is a clear necessity \citep{Peng11}. 
We refer to \citet{Harvey02}, which details the advantages of
open-source codes in hydroinformatics and claims that computer software is the lifeblood of hydroinformatics.

Some programs have adopted a strategy consisting of proposing a set of free, open-source codes to solve the shallow-water equations (e.g., GeoClaw \citep{GeoClaw11}, Gerris \citep{popinet11}, Telemac-Mascaret\footnote{\url{http://www.opentelemac.org/}} and Dassflow\footnote{University of Toulouse, CNRS, INSA. \url{http://www.math.ups-tlse.fr/DassFlow}}).  
In an explicit effort to facilitate reproducibility in water-flow modeling and simulations,
and considering that ``they are no valuable excuses not to make the code available'' \citep{Barnes10},
the present paper describes a free software package for solving the shallow-water equations: FullSWOF\@.
The source code of FullSWOF has been made public under a license that grants the freedom to use, study, share,
and modify the code. The use of a standard version of C++ (ANSI) helps in increasing accessibility to
numerous users in education, science, and industry.

FullSWOF has been developed through a joint effort between mathematicians and hydrologists.
Using a set of analytic solutions to the shallow-water equations detailed in \citet{SWASHES13}, FullSWOF includes
a validation procedure that guarantees reproducibility to the users and non-regression to the developers. This procedure is, together with the version control system,
bug tracking, etc., part of the quality assurance of
FullSWOF\@. These procedures also facilitate contributions by third parties and make their inclusion clearly identifiable.
To make the code's use and development easier,
a graphical user interface and both a one-dimensional version and a two-dimensional version have been released. FullSWOF can also be included in third-party software, as already realized in \mbox{openLISEM}  version~1.67 \citep{Baartman13}.

Another specific feature of the implementation of FullSWOF is its modular architecture. 
Other open-source developments with a modular approach that couples
hybrid numerical methods (namely, the finite element and finite volume methods)
are presented in \citet{Kolditz08} and the references therein. 
Based on the presentation of NOAH (Newcastle Object-oriented Advanced Hydroinformatics), 
\citet{Kutija07} explained why an object-oriented paradigm provides software that has a high 
computational efficiency and that is easy to maintain and extend.

As in the GeoClaw and Gerris programs, FullSWOF makes use of finite volume
methods, but specific features make it better oriented toward applications in hydrology, hence its name: \emph{Full S}hallow \emph{W}ater equations for \emph{O}verland \emph{F}low.
FullSWOF makes use of the finite volume method,
which is preferred to the finite difference method because it ensures
both mass conservation 
and the positivity of water depths. A well-balanced scheme guarantees the preservation of steady states.
Special attention has been paid to specific hydrological features: 
transitions between wet and dry areas, small water depths, and various friction models.
The software program incorporates rainfall and infiltration and enables the direct use of digital topographic grids. Some of these features were announced at a conference held in 2012 \citep{Delestre14a}.

After presenting the physically based model (section~\ref{Sec:Model}) and its most important properties (section~\ref{Sec:Prop}),
the numerical methods are described in section~\ref{Sec:NumMeth}.
Then, we proceed with the description of the FullSWOF program in section~\ref{Sec:SoftDesc}.
Finally, in section~\ref{Sec:Valid}, comparisons are made with explicit solutions representative of a wide variety of flow conditions,
and application to three real cases --- a laboratory experiment, rain on a field plot, and a dam break ---  are reported.

\section{Model: Shallow-Water Equations}\label{Sec:Model}

\subsection{General Settings}

\begin{figure}
	\begin{center}
		\includegraphics[width = 0.5\textwidth]{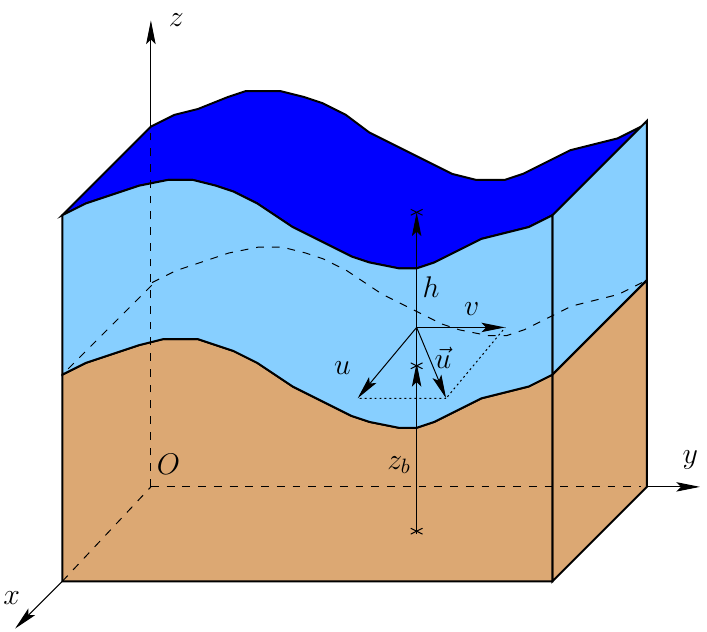}
		\caption{Notations for 2D shallow-water equations.\label{Fig:not2d}}
	\end{center}
\end{figure}

The system of shallow-water equations is a simplified model for a class of free boundary, incompressible Navier-Stokes flows
that can occur not only in rivers or channels and in the ocean (tides, tsunamis, etc.) but also in overland flow.
They are characterized by the fact that the water depth $h(t,x,y)$ [L] is small with respect to the horizontal dimensions of 
the considered domain
 (Figure~\ref{Fig:not2d} and \citet[chap. 2]{Hervouet07}).
In this context, two main hypotheses are assumed. First, the fluid
velocity is constant along the vertical direction; therefore, we can
use the horizontal components of the vertically averaged velocity $u(t,x,y)$ and $v(t,x,y)$  [L/T] instead of 
the three-component Navier-Stokes velocity vector.
Next, the pressure of the fluid is hydrostatic; therefore, after integration along the vertical direction $z$, the pressure field
is given by $p(t,x,y)=gh(t,x,y)^2/2$, where $g$ is the gravity constant [L/T$\rm^2$].

Under these assumptions, the averaged Navier-Stokes system can be rewritten as the following balance laws:
\begin{equation}\label{Eq:SaintVenant2D}
	\left\{
	\begin{array}{rcl}
		\partial_t h +\partial_x\left(hu\right)+\partial_y\left(hv\right)&=&R-I,\\[0.2cm]
		\partial_t\left(hu\right)+\partial_x\left(hu^2+\dfrac{gh^2}{2}\right)+\partial_y\left(huv\right)
  			&=& gh({S_0}_x-S_{fx}),\\[0.3cm]
		\partial_t\left(hv\right)+\partial_x\left(huv\right)+\partial_y\left(hv^2+\dfrac{gh^2}{2}\right)
			 &=& gh({S_0}_y-S_{fy}).
	\end{array}
	\right.
\end{equation}
The first equation is the exact integrated form of the incompressibility condition, and hence, it is a mass balance. 
The other two equations are momentum balances and involve forces such as gravity and friction.
In particular, the $gh^2/2$ term is the hydrostatic pressure.
Let us now describe each term, recalling their physical dimensions.
\begin{enumerate}
	\item $z_b$ is the topography [L]. Because we do not consider erosion, it is a given function of space,
		$z_b(x,y)$, and we classically denote by 
		$S_{0x}$ and $S_{0y}$ the opposites of the slopes in the $x$ and $y$ directions, respectively:
		$S_{0x}=-\partial_x z_b(x,y)$ and $S_{0y}=-\partial_y z_b(x,y)$.
	\item $R$ is the rain intensity [L/T]. It is a given function $R(t,x,y)\ge 0$.
	In the current versions of FullSWOF, we consider the rain to be uniform in space.
	\item $I$ is the infiltration rate [L/T]. This term, $I(t,x,y)\ge 0$, is defined through the coupling with an infiltration model
 such as the bi-layer Green-Ampt model (section~\ref{Sec:Infiltration}).
	\item $S_f=\left(S_{fx},S_{fy}\right)$ is the friction force, which is, in general, a nonlinear function of the
velocity and water depth (section~\ref{Sec:Friction}).
\end{enumerate}

We shall pay particular attention to the one-dimensional version of system~\eqref{Eq:SaintVenant2D} because, on the one hand, it has 
practical applications when the flow can be considered homogeneous
or when the effect of the edges can be neglected (e.g., wide channels or flood propagation in river networks).
On the other hand, its study gives a better insight into the complete two-dimensional
model from both theoretical and numerical perspectives. This one-dimensional system can
be written as
\begin{equation}\label{Eq:SaintVenant1D}
	\left\{\begin{array}{rcl}
       		\partial_t h+\partial_x (hu) &=& R-I,\\[0.2cm]
  		\partial_t (hu) +\partial_x \left(hu^2+\dfrac{gh^2}{2}\right)&=&gh(S_{0x}-S_{fx}).
      	 \end{array}\right.
\end{equation}

In both systems~\eqref{Eq:SaintVenant2D}
and~\eqref{Eq:SaintVenant1D},
the homogeneous part, that is, the left-hand side, is called the transport (or convection) operator.
This part corresponds to the flow of an ideal fluid on a flat bottom, without friction, rain or infiltration. 
In the one-dimensional setting, this is exactly the model introduced by \citet{saintvenant71}.
This operator contains several important properties of the flow; hence, in sections~\ref{Sec:Prop} and~\ref{Sec:NumMeth}, we perform
 a careful analysis of the homogeneous system in the following, considering
 both theoretical and numerical studies.

\subsection{Friction Terms}\label{Sec:Friction}

Friction terms depend on the flow velocity.
In the formul\ae\ below, $\vec{u}$ is the velocity vector $\vec{u}=(u,v)$, with $|\vec{u}| = \sqrt {u^2 + v^2}$, and $\vec{q}$ 
is the discharge $\mbox{$\vec{q}=(hu,hv)=h \vec{u}$}$.
In hydrological models, based on empirical considerations, friction
laws can be categorized into two families.
On the one hand, the Manning-Strickler friction law reads
\begin{equation}\label{Eq:Manning}
	S_f=C_f\dfrac{\vec{u}|\vec{u}|}{h^{4/3}}=C_f\dfrac{\vec{q}|\vec{q}|}{h^{10/3}},
\end{equation}
with $C_f=n^2$ or $C_f=1/K^2$, where $n$ is the Manning coefficient [L\textsuperscript{-1/3}T] and $K$ is the Strickler coefficient  [L\textsuperscript{-1/3}T\textsuperscript{-1}]
\citep{Chow59}.
On the other hand, the laws of Darcy-Weisbach and Ch\'ezy can be
written as
\begin{equation}\label{Eq:Darcy}
	S_f=C_f\dfrac{\vec{u}|\vec{u}|}{h}=C_f\dfrac{\vec{q}|\vec{q}|}{h^{3}}.
\end{equation}
Taking $C_f=f/(8g)$, with $f$ as a dimensionless coefficient, or $C_f=1/{C^2}$, $C$ [L\textsuperscript{1/2}T\textsuperscript{-1}], 
we obtain the Darcy-Weisbach and Ch\'ezy friction laws,
respectively. Readers are referred to chapter~5 of \citet{Chow59} for
details about and examples of friction laws.
Note that the friction force may depend on the space variable, especially in large domains, but this is not considered in the following.

\subsection{Infiltration Model}\label{Sec:Infiltration}

\begin{figure}[htbp]
	\subcaptionbox{Wetting front.\label{Fig:GA:wet}}[0.5\textwidth]
			{\includegraphics[width=0.5\textwidth]{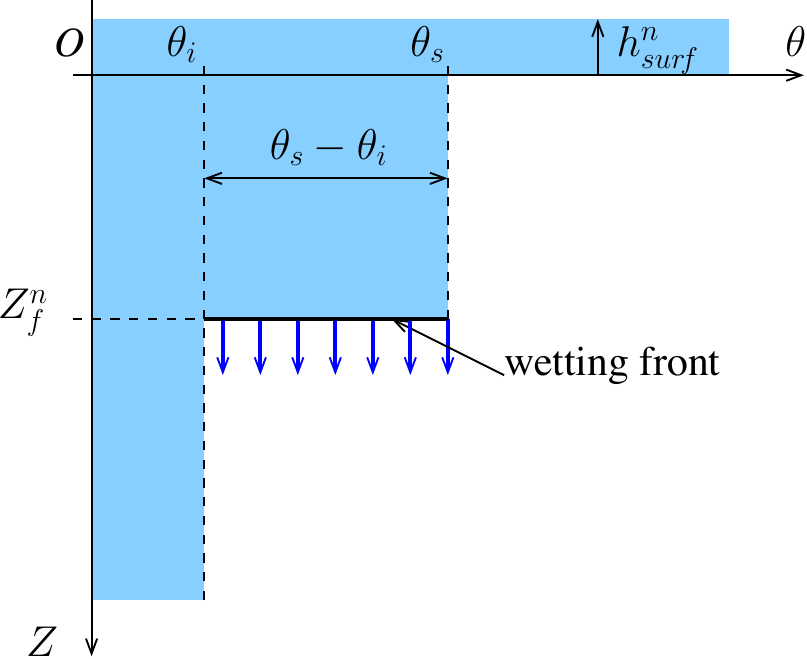} \\ \vspace{3.8em}} %
	\subcaptionbox{Geometry of the bi-layer model.\label{Fig:GA:dry}}[0.5\textwidth]
		{\includegraphics[width=0.2\textwidth]{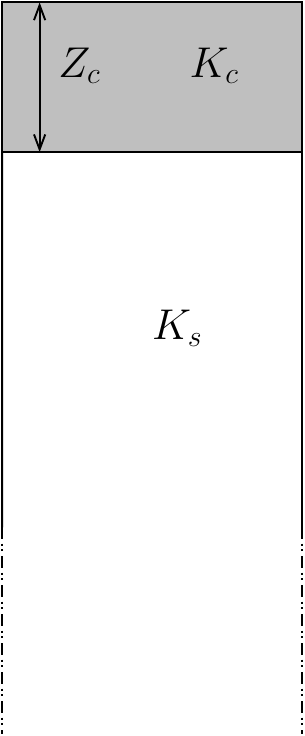}}
	\caption{Notations for the Green-Ampt infiltration model.\label{Fig:GA}}
\end{figure}

Infiltration is computed at each cell using a Green-Ampt model \citep{GreenAmpt11, Mein73}.
The main idea is to assume that the water infiltrates as an advancing wetting front (at the depth $Z_f = Z_{f}(t)$) from 
a fully saturated zone (with moisture content $\theta_s$) to another zone 
with an initial water content $\theta_i$ (Figure~\ref{Fig:GA:wet}). 

Following \citet{Esteves00}, we implemented a bi-layer Green-Ampt model \citep{Hillel70, Delestre10b} in which the 
upper layer is characterized by its thickness $Z_c$ and its hydraulic conductivity $K_c$, 
and the second layer has an infinite extension and a hydraulic conductivity $K_s$ (Figure~\ref{Fig:GA:dry}). 
Note that the `c' and `s' subscripts stand for `crust' and `soil', respectively. 
$Z_c$ can be set equal to zero for cases where a single layer model is relevant. In the following, we assume that the infiltration parameters can vary in space but are uniform in time.

At each cell, when there is water at the surface, the infiltration capacity $I_C$ [L/T] at time $t_n$ is given by
\begin{equation}
	I_{C} (t_{n}) = I_C^n = \left\{\begin{array}{l c l}
              K_s\left(1+\dfrac{h_f-h_{surf}^n}{Z_f^n}\right) & \text{ if }& Z_c = 0,\\[0.5cm]
              K_c\left(1+\dfrac{h_f-h_{surf}^n}{Z_f^n}\right) & \text{ if }& Z_f^n\leq Z_c,\\[0.5cm]
              K_e^n\left(1+\dfrac{h_f-h_{surf}^n}{Z_f^n}\right) & \text{else.}&
             \end{array}\right. 
\end{equation}
In these equalities, $K_e^n = K_{e}(t_{n})$ is the effective hydraulic conductivity at time $t_n$:
\begin{equation*}
K_e^n=\dfrac{Z_f^n}{\dfrac{Z_f^n-Z_c}{K_s}+\dfrac{Z_c}{K_c}} 
	= \dfrac{1}{\dfrac{1}{K_s}\left(1- Z_c \dfrac{\Delta \theta}{V_{in\!f}^n}\right)+ Z_c \dfrac{\Delta \theta}{V_{in\!f}^n}\dfrac{1}{K_c}}.
\end{equation*}
The value of $h_f$ depends on the soil; it is sometimes denoted by $\Psi$ in the literature and represents the suction head at the wetting front. 
The (positive) quantity $-h_{surf}^n = - h_{surf}(t_{n})$ is the water depth at the surface of the cell 
that is available for infiltration at time $t_{n}$. 
Finally, $Z_f^n = Z_{f}(t_{n})$ is the gravity force over the water column at time $t_{n}$ and can be written as 
$Z_f^n = V_{in\!f}^n/ \Delta \theta$, where 
$V_{in\!f}^n$ is the infiltrated volume at time $t_{n}$ and $\Delta \theta = \theta_s - \theta_i$. 

To avoid an infinite infiltration rate initially (when the infiltrated volume is still equal to zero),
we add a threshold to obtain the infiltration rate $I^n = \min(I_C^n, i_{max})$.
Because the infiltrated volume cannot exceed the water depth, the volume is updated as follows:
\begin{equation*}
	V_{in\!f}^{n+1} = V_{in\!f}^n + \min (-h_{surf}^n, I^n \times \Delta t)  
\end{equation*}
Finally, the water depth is updated.

\section{Properties}\label{Sec:Prop}

In this section, we recall several mathematical properties of the shallow-water model 
for both the 1D and 2D cases. 

\subsection{One-Dimensional Model}\label{Sec:1D}

To emphasize the mathematical properties of the shallow-water model, we first rewrite the one-dimensional homogeneous equations using vectors:
\begin{equation}\label{Eq:Nonconservative}
	\partial_t W+\partial_x F(W) = 0,\  \mbox{ where } W=\left(\begin{array}{c}h\\hu\end{array}\right), \  
  		F(W)=\left(\begin{array}{c}hu\\hu^2+\dfrac{gh^2}{2}\end{array}\right),
\end{equation}
with $F(W)$ being the flux of the equation. The transport is more clearly evidenced in the following non-conservative form:
\begin{equation*}\label{Eq:MatrixTransp}
	\partial_t W+A(W)\partial_x W = 0, \quad A(W)=F'(W)=\left(\begin{array}{cc} 0 & 1 \\ -u^2+gh & 2u\end{array}\right),
\end{equation*}
where $A(W)$ is the matrix of transport coefficients.
More precisely, when $h>0$, the matrix $A(W)$ turns out to be diagonalizable, with eigenvalues
$\lambda_1(W) = u-\sqrt{gh} \,< \, u+\sqrt{gh} = \lambda_2(W)$.
This important property of having two real and distinct eigenvalues is
called strict hyperbolicity (e.g., see \citet{Godlewski96} and
references therein for more details). The eigenvalues are the
velocities of the surface waves of the fluid, 
which are the fundamental characteristics of the flow. Note that the eigenvalues coincide if $h=0$, 
that is, for dry zones. In this case, the system is no longer hyperbolic,
which induces difficulties at both the theoretical and numerical levels. 

From these formul\ae\, we recover a useful classification of flows based on the relative values of the velocities of the fluid, $u$, and of the 
waves, $\sqrt{gh}$. If  $|u| < \sqrt{gh}$, the characteristic velocities $u-\sqrt{gh}$ and $u+\sqrt{gh}$ have opposite signs, and information propagates upward as
well as downward; the flow is said to be subcritical or fluvial. In contrast, when $|u| > \sqrt{gh}$, all the information propagates downward, and the flow is said to be supercritical or torrential. 

This classification has consequences for the numerical scheme. Because we have two unknowns $h$ and $u$ 
(or, equivalently, $h$ and $q=hu$), a subcritical flow is determined by one upstream value and one downstream value, 
whereas a supercritical flow is completely determined by the two upstream values. 
Thus, for numerical simulations, we use only one of the two variables
for a subcritical inflow/outflow boundary. 
For a supercritical inflow boundary,
we have to set both variables, and for a supercritical outflow boundary, the Neumann free-boundary conditions are considered \citep[e.g.,][]{Bristeau01}.
In this context, it is useful to be able to determine whether the flow is subcritical or supercritical. 
To this end, we can consider two quantities. 
The first quantity is the Froude number:
\begin{equation*}\label{Eq:DefFroude}
	{\rm Fr} = \dfrac{|u|}{\sqrt{gh}}.
\end{equation*}
The flow is subcritical or supercritical if ${\rm Fr}<1$ or ${\rm Fr}>1$, respectively.
A more visual criterion is obtained through the so-called critical
depth $h_c$, which can be written as
\begin{equation*}\label{Eq:DefCritHeight}
	h_c = \left(\dfrac{|q|}{\sqrt{g}}\right)^{2/3},
\end{equation*}
for a given discharge $q=hu$. The flow is subcritical or supercritical if $h>h_c$ or $h<h_c$, respectively.

\subsection{Source Terms and Equilibria}

When source terms (e.g., topography, rain or friction) are involved, other properties have to be considered,
in particular, the occurrence of steady-state (or equilibrium)
solutions (solutions that do not depend on time, i.e., $\partial_t\equiv 0$).
This amounts to some balance between the flux and source terms:
\begin{equation*}
\partial_x(hu) = R-I,\quad \partial_x\left(hu^2+\dfrac{gh^2}2\right) = gh(S_{0x}-S_{fx}).
\end{equation*}

Among these solutions, the homogeneous states are of particular interest, namely, the so-called parallel flows
\begin{equation*}
\partial_x(hu) = R-I = 0, \quad gh(S_{0x}-S_{fx}) = 0,
\end{equation*}
and lakes or puddles at rest
\begin{equation*}
u = 0, \quad h+z = \mbox{\rm const}.
\end{equation*}

These solutions are important for two reasons: on the one hand, specific numerical methods have to be designed to obtain
these solutions; on the other hand, these models furnish several explicit solutions that can be used as test cases for numerical methods 
(see section~\ref{Sec:Valid} and \citet{SWASHES13}).

\subsection{Two-Dimensional Model}

The two-dimensional shallow-water system~\eqref{Eq:SaintVenant2D} can be written under the following conservative form:
\begin{equation}
	\partial_t U+\partial_x G(U)+\partial_y H(U)=S(U,t,x,y),\label{Eq:SaintVenant2Dvectoriel}
\end{equation}
where
\begin{equation*}
	U=\begin{pmatrix}h\\hu\\hv \end{pmatrix},\quad
	G(U)=\begin{pmatrix}hu\\ hu^2+\dfrac{gh^2}{2}\\huv\end{pmatrix},\quad
	H(U)=\begin{pmatrix}hv\\huv\\ hv^2+\dfrac{gh^2}{2}\end{pmatrix}
\end{equation*}
and
\begin{equation*}
	S(U,t,x,y)=\begin{pmatrix}R-I\\gh({S_0}_x-S_{fx})\\gh({S_0}_y-S_{fy})\end{pmatrix}.
\end{equation*}

If we denote by $DG_x$ and $DH_y$
the Jacobian matrices of the fluxes, namely,
\begin{equation*}
	DG_x= \begin{pmatrix} 0  & 1  & 0\\ -u^2+gh  & 2u & 0\\-uv & v  & u\end{pmatrix}
	\mbox{\quad and\quad}
	DH_y=\begin{pmatrix} 0  & 0 & 1\\ -uv & v & u\\ - v^2+gh & 0 & 2v \end{pmatrix},
\label{Eq:matricesJacobiennes2D}
\end{equation*}
then system~\eqref{Eq:SaintVenant2Dvectoriel} reads
\begin{equation*}
	\partial_t U+DG_x\,\partial_x U+DH_y\,\partial_y U = S(U,t,x,y).\label{Eq:hyperbolique2d}
\end{equation*}

The notion of hyperbolicity is defined here, following \citet{Godlewski96}, by studying the flow along any direction. For any unit vector
$\xi = \left( \xi_x, \xi_y \right) \in \mathbb{R}^2$,
the velocity of the flow in the $\xi$ direction is, by definition, $u_{\xi}=\xi_x u+\xi_y v$, and we define a type of
directional derivative as $DF(\xi)=\xi_x DG_x+\xi_y DH_y$. Then, it can be verified that $DF(\xi)$ has three eigenvalues:
\begin{equation}\label{Eq:eigenv}
	\lambda_1(\xi)=u_{\xi}-\sqrt{gh},\quad \lambda_2(\xi)=u_{\xi} \mbox{\quad and \quad} \lambda_3(\xi)=u_{\xi}+\sqrt{gh}.
\end{equation}

Outside the dry zones, that is, when $h>0$, the three eigenvalues satisfy $\lambda_1(\xi)<\lambda_2(\xi)<\lambda_3(\xi)$; therefore,
the matrix $DF(\xi)$ is diagonalizable for all $\xi$. Specifically, the system~\eqref{Eq:SaintVenant2Dvectoriel} is strictly hyperbolic. As 
in the one-dimensional case, this property is no longer true inside the dry zones. 

\section{Numerical Methods}\label{Sec:NumMeth}

In this section, we detail the numerical methods that are shown to be well
adapted to finding solutions of the shallow-water system and that are
implemented in FullSWOF\@.
To obtain the solutions, the shallow-water system is divided into two parts: the transport (i.e., convective) operator and the source terms.
First, we perform a convective step, where the homogeneous part of the system is solved
by a finite volume strategy. This leads to a first-order-accurate
scheme, the second-order accuracy
in space being obtained by reconstruction techniques.
Second, this scheme is coupled with the source terms, and specific methods are used to consider
steady states.
In case a second-order approximation in time is requested, Heun's method is used.
These steps are described below in detail for the one-dimensional
setting and are then extended to the two-dimensional case.

\subsection{General Settings}
\begin{figure}[htbp]
	\subcaptionbox{One-dimensional time and space discretization.\label{Fig:discretisations:1d}}
		{\includegraphics[width=0.45\textwidth]{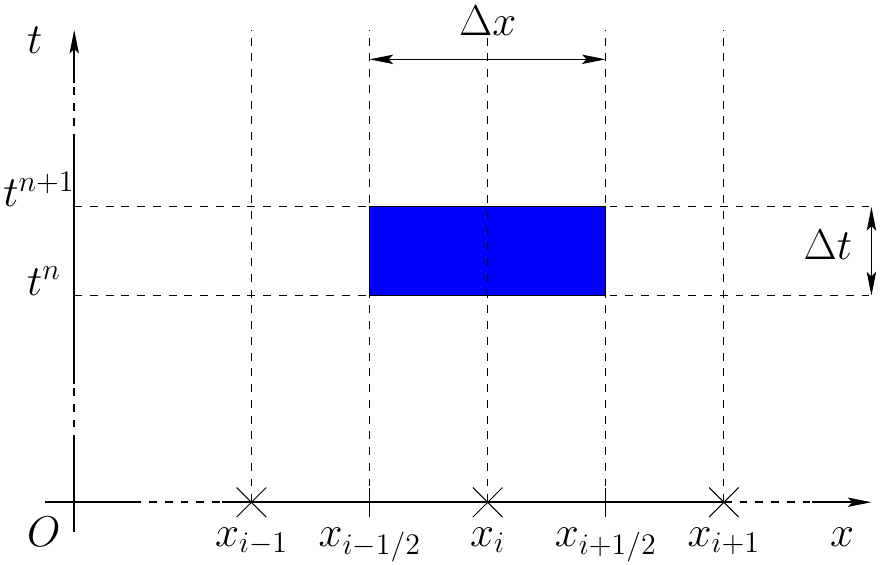}}%
	\hspace{0.10\textwidth}
	\subcaptionbox{Two-dimensional space discretization for every time step.\label{Fig:discretisations:2d}}
		{\includegraphics[width=0.45\textwidth]{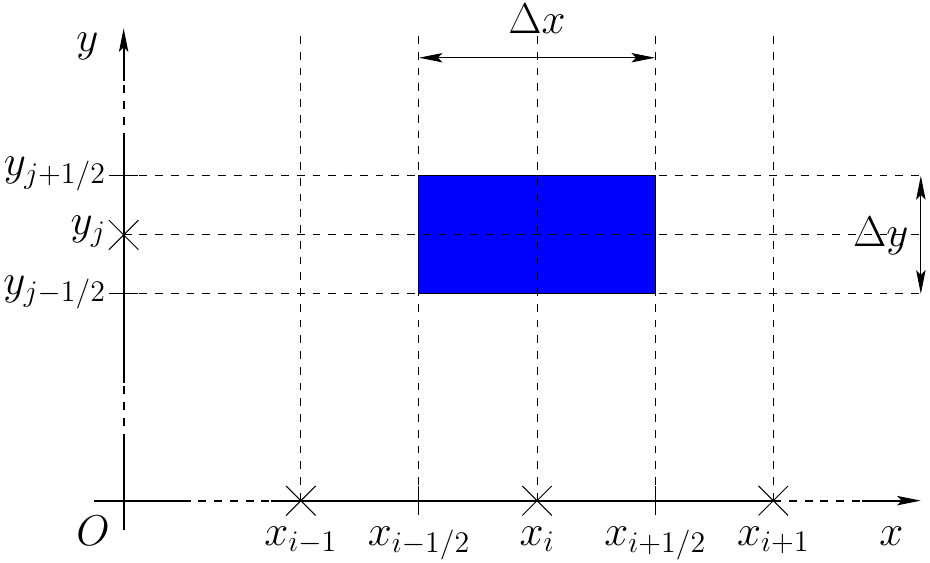}}
\caption{Discretization of time and space in FullSWOF\@.\label{Fig:discretisations}}
\end{figure}

Let us first introduce several notations that will be used throughout this section. First, consider the time discretization: 
let $t^0=0$ be the initial time; we fix a time step $\Delta t>0$, and for $n\ge 0$, we set $t^{n+1}=t^n+\Delta t$.
Next, the space discretization is defined by constant positive space
steps $\Delta x$ in one dimension and $(\Delta x,\Delta y)$ in 
two dimensions, where a rectangular mesh is assumed. 
Finally, we compute piecewise constant approximations of the vectors $W$ (in 1D) or $U$ (in 2D). 
More precisely, $W_i^n$ and $U_{ij}^n$ are constant
approximations of $W$ on $[t^n,t^{n+1}[\times]x_{i-1/2},x_{i+1/2}[$ and of $U$ on 
$[t^n,t^{n+1}[\times]x_{i-1/2},x_{i+1/2}[\times]y_{j-1/2},y_{j+1/2}[$, 
respectively (see Figure~\ref{Fig:discretisations} for details of the notations).

For each $n\ge 1$, we compute the components of
$W^n$ and $U^n$ using explicit schemes in time that take the general recursive form
\begin{equation}
W_i^{n+1}=W_i^n-\Delta t \Phi(W^n),\quad U_{ij}^{n+1}=U_{ij}^n-\Delta t \Phi(U^n).\label{Eq:generalscheme}
\end{equation}
Formula~\eqref{Eq:generalscheme} is nothing more than an Euler scheme in time, where the
function $\Phi$ describes the discretization of the space derivatives and of the source terms.

At this point, we briefly recall several fundamental properties required in an effective numerical scheme
(for detailed definitions, see \citet{Godlewski96}). 
First, the consistency error quantifies the ability of the scheme to mimic the equation.
This error is obtained by replacing each occurrence of $W_i^n$ in the formula~\eqref{Eq:generalscheme} by the exact solution $W(t^n,x_i)$. 
If this error goes to $0$ when $\Delta t$ and $\Delta x \to 0$, the scheme is consistent. 
If the error behaves as $O(\Delta t^q+\Delta x^p)$, 
the scheme is of order $q$ in time and $p$ in space. Next, stability, a notion that
can take various forms,  is a concern. Usually, it is formulated by stating that a given norm of the numerical solution at time $t^n$
is controlled by the norm at time $0$. Regarding the shallow-water equations, a particularly important stability property
is the preservation of the positivity of the water depth, especially if wet-dry transitions and thin water layers are to be simulated.

As far as precision in time is concerned, note that formula~\eqref{Eq:generalscheme} defines a first-order scheme in time,
regardless of the definition of $\Phi$. There are several ways to increase the accuracy of such schemes following standard methods to solve
ordinary differential equations. Here, we choose Heun's method (also referred to as the modified Euler method), 
which is a prediction-correction
method. In the first two formul\ae\ below, $W^*_i$ and $W_i^{**}$ are the predicted values; the last formula is the correction step: 
\begin{equation*}
W^*_i = W_i^n-\Delta t\Phi(W^n),\quad W_i^{**} = W_i^*-\Delta t\Phi(W^*), \qquad W_i^{n+1} = \dfrac{W_i^n+W^{**}_i}2.
\end{equation*}

We conclude this paragraph by emphasizing that when $\Phi$ is determined, the numerical scheme is complete.
Therefore, the remainder of this section
explains step by step how to build the function $\Phi$, which encodes
all the stability and accuracy properties of the scheme in the space variables.

\subsection{Convective Step for the One-Dimensional Model}

To solve system~\eqref{Eq:SaintVenant1D}, the first step consists of using a finite volume strategy 
to solve the homogeneous system~\eqref{Eq:Nonconservative}.
The idea is to integrate the system of equations
in each time-space cell, as depicted by the blue rectangle
in Figure~\ref{Fig:discretisations:1d}.
We obtain the following approximation formula
\begin{equation}\label{Eq:FVscheme}
W_i^{n+1} = W_i^n - \dfrac{\Delta t}{\Delta x}(F^n_{i+1/2}-F^n_{i-1/2}),
\end{equation}
where $F^n_{i+1/2}$ and $F^n_{i-1/2}$ are approximations of the flux on the edges $x_{i+1/2}$ and $x_{i-1/2}$, respectively.
At this stage, the finite volume scheme is not completely determined yet;
we need to specify a formula to compute this flux approximation at each edge.
Here, we choose to set, for all indexes $i$, $F^n_{i+1/2} = {\cal F}(W_i^n, W_{i+1}^n)$,
where ${\cal F}(W_L,W_R)$ is called the numerical flux,
a given function of the states on the left and right of the interface. 

Among the numerical fluxes proposed in the literature,
such as Rusanov, HLL, VFRoe-ncv, and kinetic (see \citet{Bouchut04} for explicit formul\ae\ and additional references),
several fluxes are available in FullSWOF\@.
In this paper, we choose to present the HLL flux (introduced in \citet{Harten83}) because
this flux offers a good compromise between simplicity, computing time and robustness, as shown
numerically in \citet{Delestre10b}. This flux can be written as
\begin{equation}\label{Eq:HLL}
{\cal F}\left(W_L,W_R \right)=
\left\{\begin{array}{ll}
        F(W_L) & \text{if}\; 0\leq c_1\\
\dfrac{c_2 F(W_L)-c_1 F(W_R)}{c_2 -c_1}+\dfrac{c_1 c_2}{c_2 -c_1}\left(W_R - W_L\right) & \text{if}\; c_1<0<c_2 \\
        F(W_R) & \text{if}\; c_2 \leq 0
       \end{array}\right.,
\end{equation}
with two parameters $c_1<c_2$
that are the approximations of the slowest and fastest wave speeds. We refer to \citet{Batten97}
 for further discussion on the wave speed estimates. In FullSWOF, we have implemented
\begin{equation*}
c_{1}={\inf\limits_{W=W_L,W_R}}({\inf\limits_{j\in\{1,2\}}}\lambda_{j}(W))\;\text{and}
\;c_{2}={\sup\limits_{W=W_L,W_R}}({\sup\limits_{j\in\{1,2\}}}\lambda_{j}(W)), 
\end{equation*}
where $\lambda_1(W)=u-\sqrt{gh}$ and $\lambda_2(W)=u+\sqrt{gh}$ are the eigenvalues of the one-dimensional model (section~\ref{Sec:1D}).
This numerical flux is said to be upwind, i.e., it mimics the wave propagation. When the flow is supercritical ($0\leq c_1$
 or $c_2 \leq 0$), all information is going from upstream to downstream; thus, the numerical flux is calculated from upstream values only.
In contrast, for subcritical flows ($c_1<0<c_2$), information is coming from both upstream and downstream;
thus, the flux is calculated using both upstream and downstream values.

As described above, except for certain numerical fluxes, the scheme is first-order accurate in space. To obtain second-order accuracy, we perform a linear reconstruction in space on the variables $W=(h,hu)$, 
thus obtaining new variables on each interface
$i+1/2$, namely, $W_{i+1/2-}$ on the left and $W_{i+1/2+}$ on the right. 
Using ${\mathcal F}(W_{i+1/2-},W_{i+1/2+})$ in the finite volume scheme instead of ${\mathcal F}(W_i^n,W_{i+1}^{n})$ 
leads to a second-order approximation in space.
Several of the existing formul\ae\ for the linear reconstruction (e.g.,
MUSCL and ENO, see \citet{Bouchut04}) 
are implemented in FullSWOF\@. 

We conclude this section with an important remark about time discretization. Explicit schemes imply a control on the time step, 
which cannot be too large for a given 
space step. More precisely, for a 3-point scheme, as defined above ($W_i^{n+1}$ depends only on $W_{i-1}^n$, $W_i^n$ and $W_{i+1}^n$),
the numerical speed of propagation is $\Delta x/\Delta t$. 
To avoid any loss of information, this velocity has to be larger than
any possible physical velocity:
\begin{equation*}
C\dfrac{\Delta x}{\Delta t} \ge \sup_{(t,x)}\left\{|u(t,x)|+\sqrt{gh(t,x)}\right\}, 
\end{equation*}
where the supremum is taken over the whole time-space domain of interest and $C$ is a parameter depending on the dimension and on the order of
the schemes that are considered (in 1D, at the first order, $C=1$, and, at the second order, $C=0.5$; in 2D, at the first order, $C=0.5$, 
and, at the second order, $C=0.25$).
When this condition is violated, the scheme becomes unstable, and oscillations appear.
This theoretical formulation of the limitation, known as the
CFL condition (Courant, Friedrichs, Lewy) \citep{Godlewski96}, is not useful because,
due to the nonlinearity, the right-hand side is difficult to estimate. Moreover,
even if such an estimate is obtained, it can lead to an underestimated fixed time step
when there are no large variations in $h$ and $u$.
Thus, the computation of the time step could be replaced by
a sequence of variable time steps $\Delta t^n$ according to the following rule:
\begin{equation}\label{Eq:timestep}
\Delta t^{n} =  \frac{C \Delta x}{\underset{i}\sup\{|u_i^{n-1}|+\sqrt{gh_i^{n-1}}\}}.
\end{equation}
However, when water depths converge toward zero, the time steps will tend to become infinite, which is not reasonable. 
Consequently, we add an upper bound to the computation, which leads to the following formulation (implemented in FullSWOF\_2D):
\begin{equation*}
\Delta t^{n} =C \min\left( \Delta x, \frac{\Delta x}{\underset{i}\sup\{|u_i^{n-1}|+\sqrt{gh_i^{n-1}}\}}\right).
\end{equation*}

\subsection{Source Terms}

When source terms are involved, specific methods have to be introduced to obtain equilibrium states. 
Numerical schemes that preserve stationary equilibria
are called equilibrium schemes, or well-balanced schemes \citep{Bermudez94, Greenberg96}. 
The first class of methods consists of splitting-type methods, where
the transport equation and source term are solved independently. Another strategy consists of
applying a specific reconstruction technique, thus modifying the flux
computation. Apparently, the latter method is well adapted 
to the topography term but not quite to the friction term.
Hence, in FullSWOF, the following two approaches are used: the topography is solved using a well-balanced scheme (section~\ref{Sec:Hydro}), and 
the friction is solved with a semi-implicit method (section~\ref{Sec:Frictionnum}).

\subsubsection{The Hydrostatic Reconstruction}\label{Sec:Hydro}

The hydrostatic reconstruction \citep{Audusse04c,Bouchut04}, as its name suggests, reconstructs new variables to be used in the 
numerical flux. It is designed to obtain a well-balanced scheme in the sense that it preserves, at least, the steady state at rest
(i.e., the hydrostatic equilibrium) as well as the positivity of the water depth.
The hydrostatic reconstruction procedure is first applied 
to the original variables
$W_i^n$, thus leading to a first-order scheme.

Here, we detail how to obtain a second-order version, which is used in FullSWOF\@.
The second-order scheme consists of first performing the linear reconstruction (chosen for the convective step)
not only to $W=(h,hu)$, as before, but also to $h+z_b$.
Then, the hydrostatic reconstruction is applied to these modified variables.
Note that this strategy introduces an artificial time
dependence on the topography, which is in some sense reconstructed as well. This is 
mandatory for addressing  equilibrium solutions and for preserving the positivity of $h$ \citep{Audusse04c}.
The linear reconstruction gives the
values $(h_{i+1/2-},z_{i+1/2-},u_{i+1/2-})$ on the left of the interface $i+1/2$ and $(h_{i+1/2+},z_{i+1/2+},u_{i+1/2+})$
on the right of the interface.
The final formula for the hydrostatic reconstruction can now be written as
\begin{equation}
 \left\{\begin{array}{l}
         h_{i+1/2L}=\max \left(h_{i+1/2-}+z_{i+1/2-}-\max\left(z_{i+1/2-},z_{i+1/2+}\right),0\right),\\
         W_{i+1/2L}=\begin{pmatrix}h_{i+1/2L}\\h_{i+1/2L}u_{i+1/2-}\end{pmatrix},\\
         h_{i+1/2R}=\max \left(h_{i+1/2+}+z_{i+1/2+}-\max\left(z_{i+1/2-},z_{i+1/2+}\right),0\right),\\
         W_{i+1/2R}=\begin{pmatrix}h_{i+1/2R}\\h_{i+1/2R}u_{i+1/2+}\end{pmatrix}.
        \end{array}\right.\label{Eq:hydrostatic-reconstruction}
\end{equation}
A given space discretization may exhibit abnormal behaviors for some combinations of slopes and water depths 
\citep{CRAShydro}. Particularly obvious for the first-order scheme and
on a coarse mesh, they disappear when refining the mesh and are hardly
noticeable at second order.

When applying the hydrostatic reconstruction,
formula~\eqref{Eq:FVscheme} has to be modified to preserve the
consistency of the scheme. This can be written as
$W_i^{n+1}=W_i^n-\Delta t \Phi(W^n)$, with
\begin{equation*}
\Phi (W^n)= \dfrac{1}{\Delta x} \left(F^n_{i+1/2L}-F^n_{i-1/2R}-Fc^n_i \right),
\end{equation*}
where $F^n_{i+1/2L}$ and $F^n_{i-1/2R}$ are given by
\begin{equation*}\begin{cases}
F^n_{i+1/2L} = {\cal F}(W^n_{i+1/2L},W^n_{i+1/2R}) + 
	\begin{pmatrix}0\\ \dfrac{g}{2}\left(\left(h^n_{i+1/2-}\right)^2-\left(h^n_{i+1/2L}\right)^2\right)\end{pmatrix},\\
F^n_{i-1/2R} = {\cal F}(W^n_{i-1/2L},W^n_{i-1/2R}) +
	\begin{pmatrix}0\\ \dfrac{g}{2}\left(\left(h^n_{i-1/2+}\right)^2-\left(h^n_{i-1/2R}\right)^2\right)\end{pmatrix},
	\end{cases}
\end{equation*}
where $W^n_{i+1/2L}$ and $W^n_{i+1/2R}$ are computed using formula~\eqref{Eq:hydrostatic-reconstruction}.
The additional centered term $Fc^n_i$ is determined to preserve
consistency and to ensure a well-balanced scheme~\citep{Audusse04c}:
\begin{equation*}
Fc_i^n=\begin{pmatrix} 0\\ -\dfrac{g}{2}\left(h^n_{i-1/2+}+h^n_{i+1/2-}\right)\left(z^n_{i+1/2-}-z^n_{i-1/2+}\right) \end{pmatrix}.
\end{equation*}

\subsubsection{Friction}\label{Sec:Frictionnum}
Friction may be addressed in a well-balanced scheme by first introducing the
friction in the topography term and by then applying the hydrostatic reconstruction (section~\ref{Sec:Hydro}); this approach is therefore 
named ``apparent topography''. 
This method has been used, for example, to solve a shallow-water system with a Coriolis force \citep{Bouchut04} 
and with Coulomb friction \citep{Bouchut04, Mangeney07}.
The main idea is to use the modified topography $z_{app}$ defined by $z_{app}=z_{b}-b$, with $\partial_x b = S_{fx}$. 
The use of this class of methods gives rise to a well-balanced scheme
for friction as well as for the topography and thus computes 
equilibrium states correctly. However, the solutions are not
completely satisfactory for transitory situations, as noted in \citet{Delestre09} and \citet{DelestreJames09}; a spurious peak appears at the wet-dry
 front before equilibrium is reached.
Therefore, we utilize splitting methods. The explicit discretization, despite its
simplicity, is not relevant for the type of problems we are interested
in because this discretization leads
to instabilities and overestimations of the velocity at wet-dry interfaces \citep{Paquier95}.
However, the fully implicit method has a high 
computational cost and cannot be generalized to the two-dimensional system due to the nonlinear form of the friction law.
Another possibility consists of using a Strang time splitting to
obtain the second order in time calculation; however, this leads to a more complicated algorithm, and moreover, there is no 
significant gain in accuracy compared to first-order methods \citep{Liang09b}.
Finally, the best compromise between
accuracy, stability and computational complexity lies in semi-implicit methods
\citep{Fiedler00, Bristeau01, Liang09b}. We choose the semi-implicit treatment proposed in \citet{Bristeau01} 
not only because it preserves steady states at rest but also for its stability.
At first order in time, after a convective step $W^*_i = W_i^n-\Delta t\Phi(W^n)$, the value $W_i^{n+1}$ is given by
\begin{equation*}
	W_i^{n+1} = \begin{pmatrix} h_i^{n+1}\\q_i^{n+1}\end{pmatrix} = \begin{pmatrix}
		h_i^*\\
		{q}^{*}_{i}\left(1 + g n^2\Delta t \dfrac{| {q}_{i}^{n} |}{h^n_{i} \left(h^{n+1}_{i}\right)^{4/3}}\right)^{-1}
	\end{pmatrix}
\end{equation*}
for the Manning friction law~\eqref{Eq:Manning} and by
\begin{equation*}
	W_i^{n+1} = \begin{pmatrix} h_i^{n+1}\\q_i^{n+1}\end{pmatrix} = \begin{pmatrix}
		h_i^*\\
		{q}^{*}_{i}\left(1 + \Delta t \dfrac{f}{8} \dfrac{| {q}_{i}^{n} |}{h^n_{i}h^{n+1}_{i}}\right)^{-1}
	\end{pmatrix}
\end{equation*}
for the Darcy-Weisbach friction law~\eqref{Eq:Darcy}. 
Note the simplicity of the method, which gives an explicit value for $W_i^{n+1}$.

\subsubsection{Rain and Infiltration}

Unlike the friction and topography source terms, rain and infiltration produce no particular numerical difficulties, such as steady-state or stability preservation. Moreover, addressing infiltration implicitly would complicate the integration of other infiltration models (such as Richard's or Darcy's models) into FullSWOF\@. For these reasons, we have chosen to treat the rain and infiltration terms explicitly.

\subsection{Convective Step for the Two-Dimensional Model}

In general, the 2D system can be treated in the same manner as the 1D
system because the calculations are performed for each interface of each cell (and thus do not depend on the cell geometry). Digital elevation models (DEMs) (i.e., digital topographic maps) are mainly represented as structured grids.
Some DEMs can be represented as vector-based triangular networks (TINs).
Because TINs can easily be converted into DEMs, we chose to make developments dedicated to structured meshes (Figure~\ref{Fig:discretisations:2d}).
To obtain a 2D version, we perform the linear reconstruction of $U$ and $z_b$ and the hydrostatic reconstruction as in one dimension (equation~\eqref{Eq:hydrostatic-reconstruction}).
We obtain $U_{\bullet L, \bullet}$ and $U_{\bullet R, \bullet}$ along
the $x$ direction,
and we obtain $U_{\bullet , \bullet L}$ and $U_{\bullet , \bullet R}$
along the $y$ direction.

Then, the two-dimensional finite volume scheme reads
\begin{equation*}
	U^*_{i,j} = U^n_{i,j} - \dfrac{\Delta t}{\Delta x} \left(G^n_{i+1/2L,j}-G^n_{i-1/2R,j}-Gc^n_{i,j} \right)
		- \dfrac{\Delta t}{\Delta y} \left(H^n_{i,j+1/2L}-H^n_{i,j-1/2R}-Hc^n_{i,j} \right)
\end{equation*}
with
\begin{equation*}
	\begin{cases}
		G^n_{i+1/2L,j}={\cal G}\left(U^n_{i+1/2L,j},U^n_{i+1/2R,j}\right)+S^n_{i+1/2L,j},\\
		G^n_{i-1/2R,j}={\cal G}\left(U^n_{i-1/2L,j},U^n_{i-1/2R,j}\right)+S^n_{i-1/2R,j},\\
		H^n_{i,j+1/2L}={\cal H}\left(U^n_{i,j+1/2L},U^n_{i,j+1/2R}\right)+S^n_{i,j+1/2L},\\
		H^n_{i,j-1/2R}={\cal H}\left(U^n_{i,j-1/2L},U^n_{i,j-1/2R}\right)+S^n_{i,j-1/2R},
	\end{cases}
\end{equation*}
where
\begin{align*}
	S^n_{i+1/2L,j} = \begin{pmatrix}0\\ \dfrac{g}{2}\left(\left(h^n_{i+1/2-,j}\right)^2-\left(h^n_{i+1/2L,j}\right)^2\right)\\0\end{pmatrix},\\ 
	S^n_{i-1/2R,j} = \begin{pmatrix}0\\ \dfrac{g}{2}\left(\left(h^n_{i-1/2+,j}\right)^2-\left(h^n_{i-1/2R,j}\right)^2\right)\\0\end{pmatrix},\\
	S^n_{i,j+1/2L} = \begin{pmatrix}0\\ 0\\ \dfrac{g}{2}\left(\left(h^n_{i,j+1/2-}\right)^2-\left(h^n_{i,j+1/2L}\right)^2\right)\end{pmatrix},\\
	S^n_{i,j-1/2R} = \begin{pmatrix}0\\ 0\\ \dfrac{g}{2}\left(\left(h^n_{i,j-1/2+}\right)^2-\left(h^n_{i,j-1/2R}\right)^2\right)\end{pmatrix}.
\end{align*}
${\cal G}$ and ${\cal H}$ are the numerical fluxes used to obtain the solution of the homogeneous system. 
The first two components of ${\cal G}$, ${\cal G}_1$ and ${\cal G}_2$, as well as the first and third components of ${\cal H}$, 
${\cal H}_1$ and ${\cal H}_3$, are computed as the components of ${\cal F}$.
In the HLL Riemann solver (equation~\eqref{Eq:HLL}), being a two-wave
numerical flux solver, this assumption is correct only for hyperbolic systems of two equations, such as the one-dimensional shallow-water equations.
In two spatial dimensions, there are three equations and hence three
eigenvalues (equation~\eqref{Eq:eigenv}). In this case, the HLL solver
described in equation~\eqref{Eq:HLL} is not sufficiently precise
because it only involves two waves.  To address this issue, Toro proposed an extension of HLL called the HLLC solver \citep{Toro94}. The solver we propose in FullSWOF is inspired by this reference and defines ${\cal G}_3$ and ${\cal H}_2$  as
\begin{equation*}\begin{array}{rl}
 	&{\cal G}_3\left(U^n_{i+1/2L,j},U^n_{i+1/2R,j}\right) =\\\\
	&\qquad\begin{cases}
			v^n_{i+1/2L,j}{\cal G}_1\left(U^n_{i+1/2L,j},U^n_{i+1/2R,j}\right) \quad  \mbox{if } u^n_{i+1/2L,j}+u^n_{i+1/2R,j}>0,\\
			v^n_{i+1/2R,j}{\cal G}_1\left(U^n_{i+1/2L,j},U^n_{i+1/2R,j}\right) \quad  \mbox{if } u^n_{i+1/2L,j}+u^n_{i+1/2R,j}\leq0,
	\end{cases}\end{array}
\end{equation*}
\begin{equation*}\begin{array}{rl}
	&{\cal H}_2\left(U^n_{i,j+1/2L},U^n_{i,j+1/2R}\right) =\\\\
	&\qquad\begin{cases}
                   	u^n_{i,j+1/2L}  {\cal H}_1\left(U^n_{i,j+1/2L},U^n_{i,j+1/2R}\right) \quad  \mbox{if } v^n_{i, j+1/2L}+v^n_{i, j+1/2R,}>0,\\
            		u^n_{i,j+1/2R}  {\cal H}_1\left(U^n_{i,j+1/2L},U^n_{i,j+1/2R}\right) \quad  \mbox{if } v^n_{i, j+1/2L}+v^n_{i, j+1/2R,}\leq 0.
        \end{cases}\end{array}
\end{equation*}

\section{Description of the FullSWOF Software Package}\label{Sec:SoftDesc}

FullSWOF stands for ``Full Shallow-Water equations for Overland Flow''.
The names FullSWOF\_1D and FullSWOF\_2D define the one-dimensional
and two-dimensional versions, respectively. 

Having been designed to encourage research reproducibility, 
the source codes (in C++) for FullSWOF are available and can be downloaded from the websites \url{https://sourcesup.renater.fr/projects/fullswof-1d/}
for the one-dimensional version and from \url{https://sourcesup.renater.fr/projects/fullswof-2d/} for the two-dimensional version.
The common optional graphical interface is named FullSWOF\_UI (``UI'' for User Interface),
was developed in Java, and is available for download at \url{https://sourcesup.renater.fr/projects/fullswof-ui/}.
Each piece of software is distributed under the CeCILL-V2 (GPL-compatible) free software
license, which allows use of the software package without any limitation as to its field of application.
For more details, we refer to the documentation on the websites.

The structure of the source code is designed to make future development easy, especially for new developers;
for example, a new friction law can easily 
be added to the libfriction library by creating a new friction file. 
The documentation for programmers is included directly in the C++ code using doxygen specific comments \citep{doxygen}. This leads to the automatic extraction of the programmer's manual and simplifies the documenting task.
Due to the hosting on a software forge, FullSWOF obtains the benefits of systems of version control,
bug tracking, package release, etc. 

The one-dimensional program is primarily designed, of course, to perform 1D flow simulations.
It is also used as a development tool
to test new numerical methods (e.g., fluxes) and to introduce new features into the models
(such as new friction and infiltration laws). 
The new code can then be later integrated into the 
two-dimensional program; because the 2D mesh is a structured mesh, it is easy to adapt the code
from 1D to 2D. In addition, a parallelized version (MPI) of FullSWOF\_2D has been developed to run large test cases (such as the Malpasset 
dam break, section~\ref{Sec:Valid}).
Currently considered in an early stage of development, this version is
available only in a dedicated branch of the version control system.

The FullSWOF software is designed to use point-wise-defined topographies
and/or initial conditions.
Currently, several topographies and initialization values are hard coded
(e.g., a parabolic topography and the wet dam-break initial data). 
For other cases, input text files can be read. 
The friction law can be chosen between the Manning and Darcy-Weisbach laws.
Several boundary conditions are available (e.g., wall, 
Neumann, periodic, and imposed depth). 
The rain is constant in space but can vary in time.
Several numerical methods are implemented for the flux,
the linear reconstruction and the order of the scheme. 
Based on a study of numerical methods for overland flows \citep{Delestre10b},
default values have been selected for
an HLL-type method for the flux, the MUSCL formula for the reconstruction and a second-order scheme. 
We refer to the documentation for more details.

Some classical benchmarks and analytic solutions from the literature are defined in FullSWOF\@.
They have been chosen among those gathered in SWASHES, a compilation of shallow-water analytic solutions for hydraulic and environmental studies \citep{SWASHES13}.
They are used to validate each new version of the code before its
release, assuring the quality of the software. This validation is partially automated due to a dedicated script and is useful to both users
(who can verify that their results do not differ from reference results --- and do not depend on
their compiler, operating system or hardware) and developers
(who can verify that their changes do not cause a regression in the result quality).
Additional test cases can be added to validate new conditions introduced in the code.

\section{Numerical Illustrations}\label{Sec:Valid}

In this section, we illustrate some results from FullSWOF using a few classical test cases, using analytic solutions available
in the literature, and  using three real datasets. This is by no means an analysis of the performance of FullSWOF\@. 
We merely attempt to demonstrate the ability of
FullSWOF to simulate a wide range of flow conditions. In particular, real datasets (section~\ref{Sec:Realdata}) are used
 without any filtering of the data (no smoothing, etc.). 
Note that the illustrated test cases and analytic solutions are part
of the benchmark set included in the FullSWOF\_1D or FullSWOF\_2D codes. In the following, FullSWOF\_1D version~1.01.00 (2013-05-17) and FullSWOF\_2D version~1.04.08 (2013-11-07) have been used.

\subsection{Classical Test Cases}\label{Sec:Classical}
We begin by running FullSWOF on three analytic cases: one steady-state and two classical transitory solutions. These examples
correspond to the original shallow-water system with a spatially varying topography; there is no friction, no rain and no infiltration.

\subsubsection{Lake at Rest with an Emerged Bump}\label{Sec:Bump}
This one-dimensional steady-state solution has been developed in
\citet{Delestre10b} and in \citet[§3.1.2]{SWASHES13} 
as a test case for the preservation of steady states and the boundary condition treatment.
This solution is based on the topography provided by \citet{Goutal97}. The initial condition satisfies the hydrostatic equilibrium
\begin{equation}
h+z= \mbox{\rm const}\;\text{and}\;q=0\;\si{m^2.s^{-1}}.\label{Eq:hydeq}
\end{equation}
The domain length is set to $L = 25\;\si{m}$, with a topography given by
\begin{equation*}
z(x)=\left\{\begin{array}{ll}
0.2-0.05(x-10)^2&\text{if}\;8\;\si{m}<x<12\;\si{m},\\
0&\text{else}. \end{array}\right.
\end{equation*}

\begin{figure}[htbp]
	\subcaptionbox{Water depths and topography.\label{Fig:BumpRestEmerged:h}}
		{\includegraphics[width=0.48\textwidth]{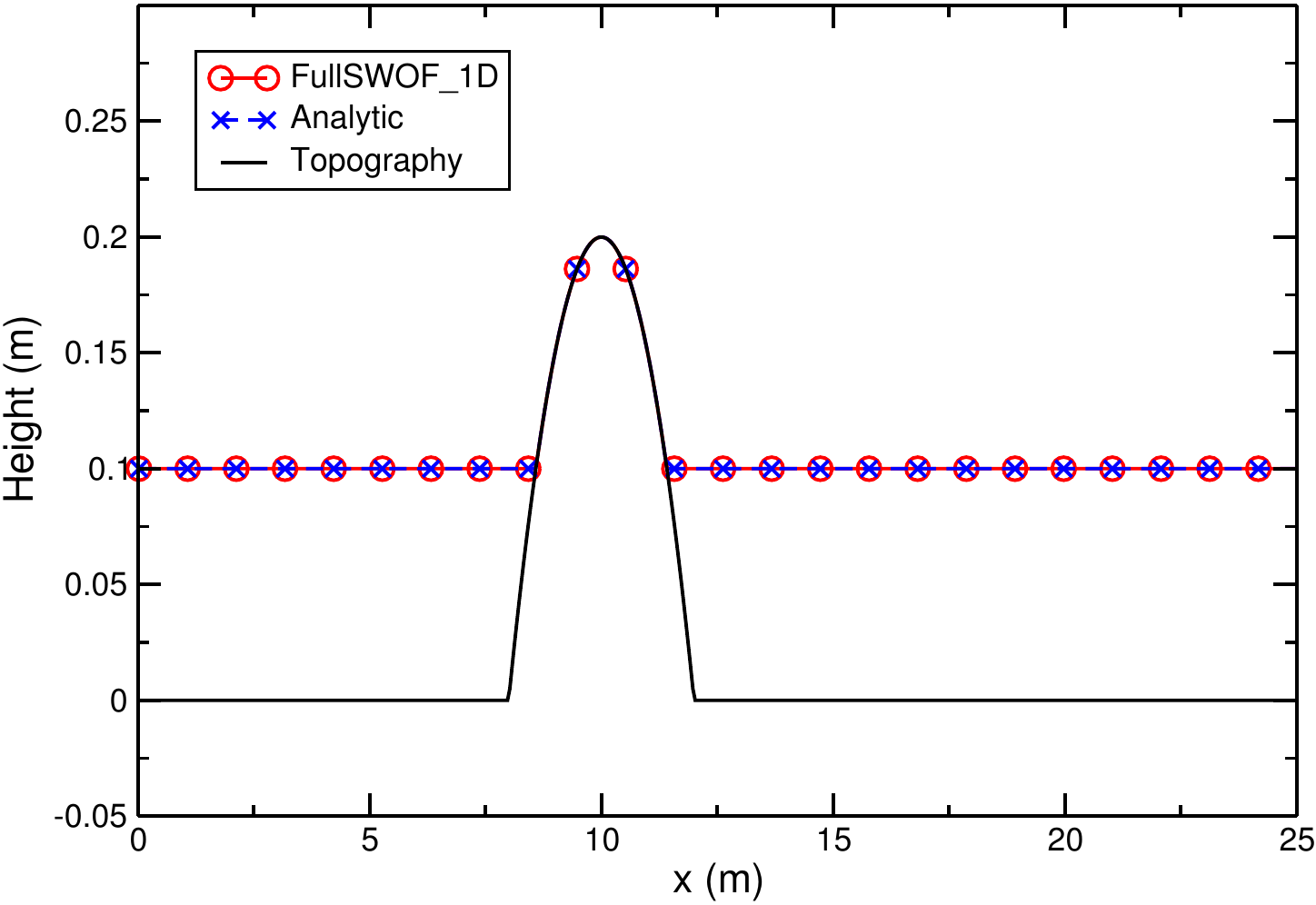}}
	\hfill
	\subcaptionbox{Water flow rates.\label{Fig:BumpRestEmerged:q}}
		{\includegraphics[width=0.48\textwidth]{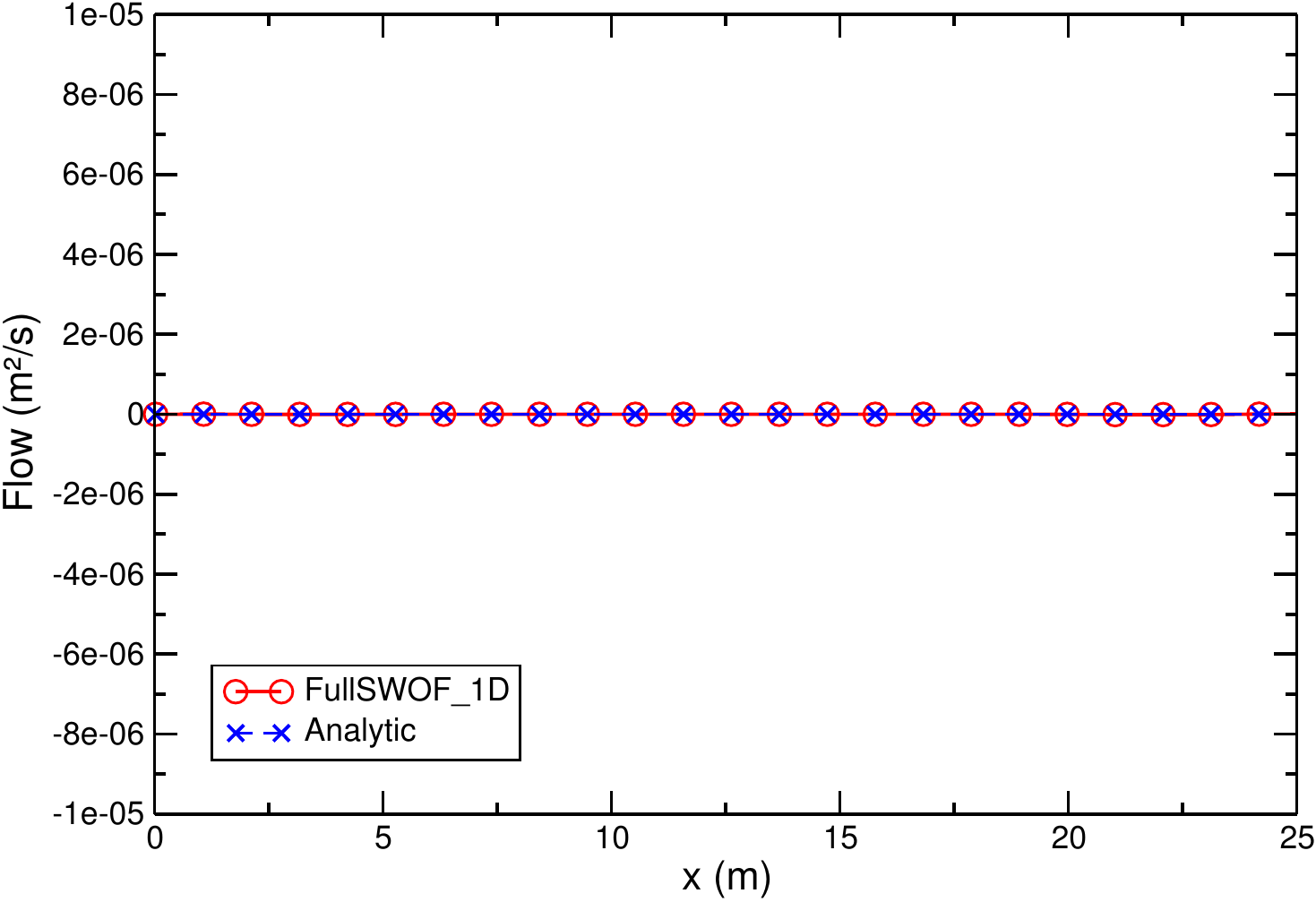}}
	\caption{Lake at rest with an emerged bump:
		comparison of the analytic solution with the FullSWOF\_1D results.
		Simulation with 500~cells at T~= 100~\si{s}.\label{Fig:BumpRestEmerged}}
\end{figure}

The maximum water depth is smaller than the amplitude of the topography to simulate a lake at rest with an emerged bump
 (Figure~\ref{Fig:BumpRestEmerged:h}).   
In such a configuration, starting from the steady state, the velocity must remain null, and the water surface should stay 
flat. This is the exact behavior simulated by FullSWOF\_1D, showing the interest in a well-balanced scheme (Figure~\ref{Fig:BumpRestEmerged}).

\subsubsection{Dam Break on a Dry Domain without Friction}

We turn now to a transitory one-dimensional case, namely, the analytic solution of a dam break on a dry domain without friction on a flat and horizontal topography \citep[§4.1.2]{SWASHES13}. This case is known as
Ritter's solution \citep{Ritter92, Hervouet07}. The difficulties here are (a)~the existence of a shock and (b)~a wet-dry transition.
This case also tests whether the scheme preserves the positivity of
the water depth because this property is usually violated near the wetting front. 

The initial condition for this configuration is the following Riemann problem
\begin{equation*}
h(x)= \begin{cases}
h_{l} = 0.005\;\si{m}& \mbox{for\quad} 0\;\si{m}\leq x \leq 5\;\si{m},\\
h_{r} = 0\;\si{m} & \mbox{for\quad} 5\;\si{m} <x\leq 10\;\si{m}, \end{cases}
\end{equation*}
with $u(x)= 0\;\si{m.s^{-1}}$.

\begin{figure}[htbp]
	\subcaptionbox{Water depths and topography.\label{Fig:DryDamBreak:h}}
		{\includegraphics[width=0.48\textwidth]{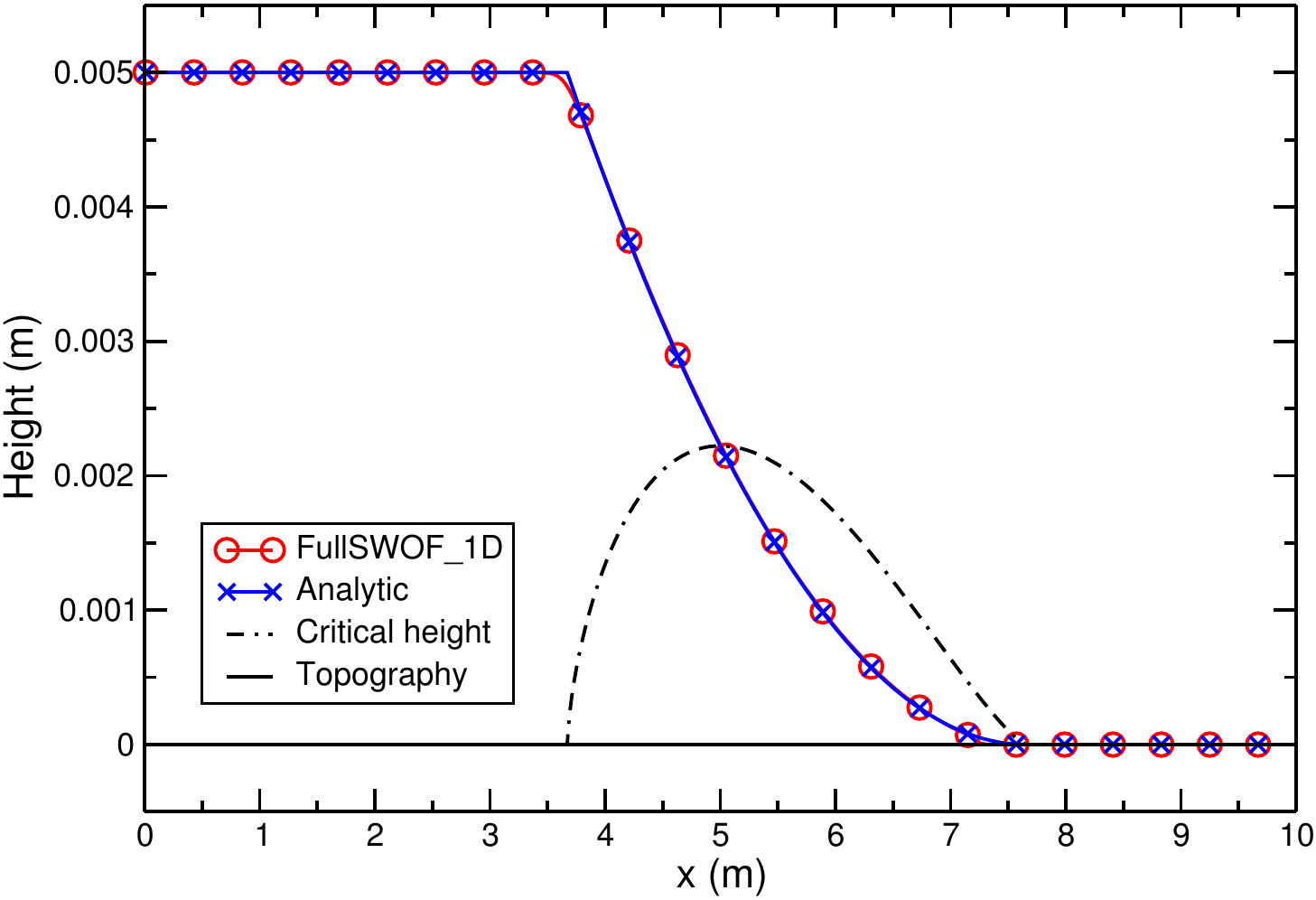}}
	\hfill
	\subcaptionbox{Water flow rates.\label{Fig:DryDamBreak:q}}
		{\includegraphics[width=0.48\textwidth]{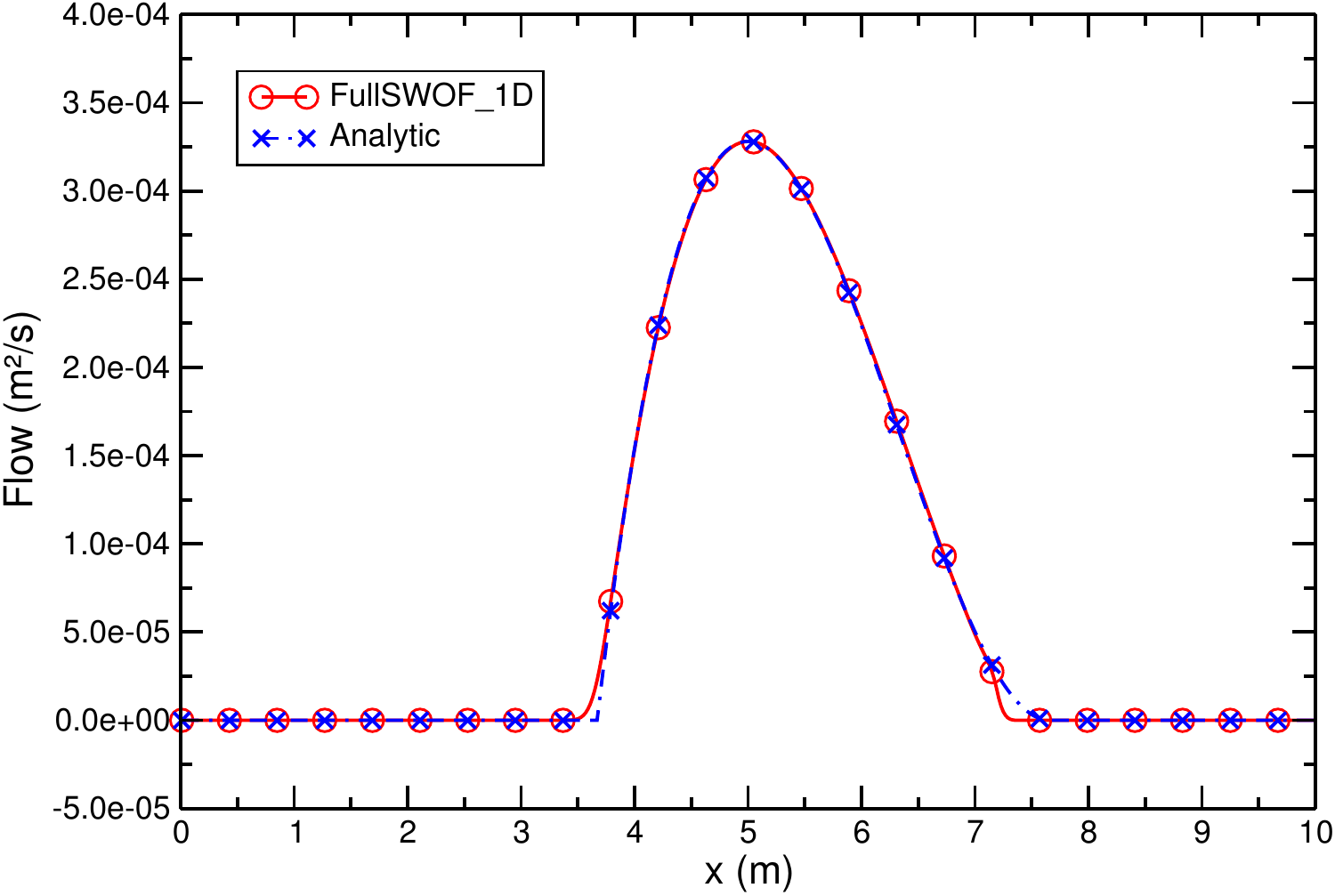}}
	\caption{Dam break on a dry domain without friction:
		comparison of the analytic solution with the FullSWOF\_1D results.
		Simulation with 500~cells at T~= 6~\si{s}.\label{Fig:DryDamBreak}}
\end{figure}

Initially, the free surface exhibits the following structure: starting from upstream, there
is a constant water depth at rest (with $h = 0.005\;\si{m}$) connected by a parabola to a dry zone downstream
(with $h = 0\;\si{m}$).
The left extremity of the parabola moves upstream, while its right end slides downstream. Figure~\ref{Fig:DryDamBreak} displays this solution at time $t=6\;\si{s}$.

There is an overall agreement between the FullSWOF\_1D results and the analytic solution;
the code is able to represent the shock, locate and correctly treat the wet-dry transition, and preserve the positivity of the water depth.
However, differences are observed at the two connections between the constant states and the parabola. 
These differences become smaller when the space step is decreased (results not shown).

\subsubsection{Planar Surface Rotating in a Paraboloid}
\begin{figure}[htbp]
	\subcaptionbox{Water depths and topography.\label{Fig:Thacker:h}}
		{\includegraphics[width=0.48\textwidth]{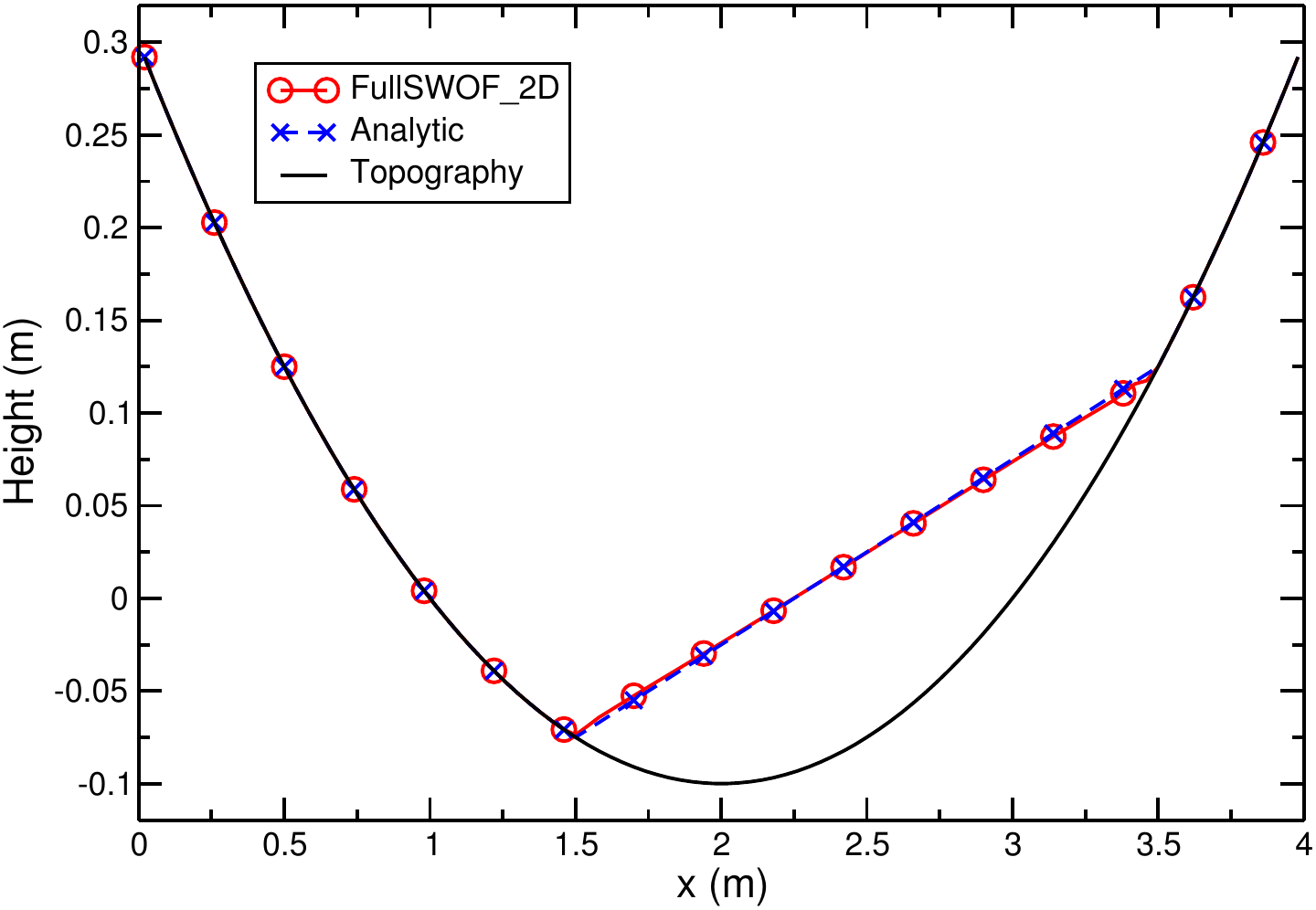}}
	\hfill
	\subcaptionbox{Water flow rates.\label{Fig:Thacker:q}}
		{\includegraphics[width=0.48\textwidth]{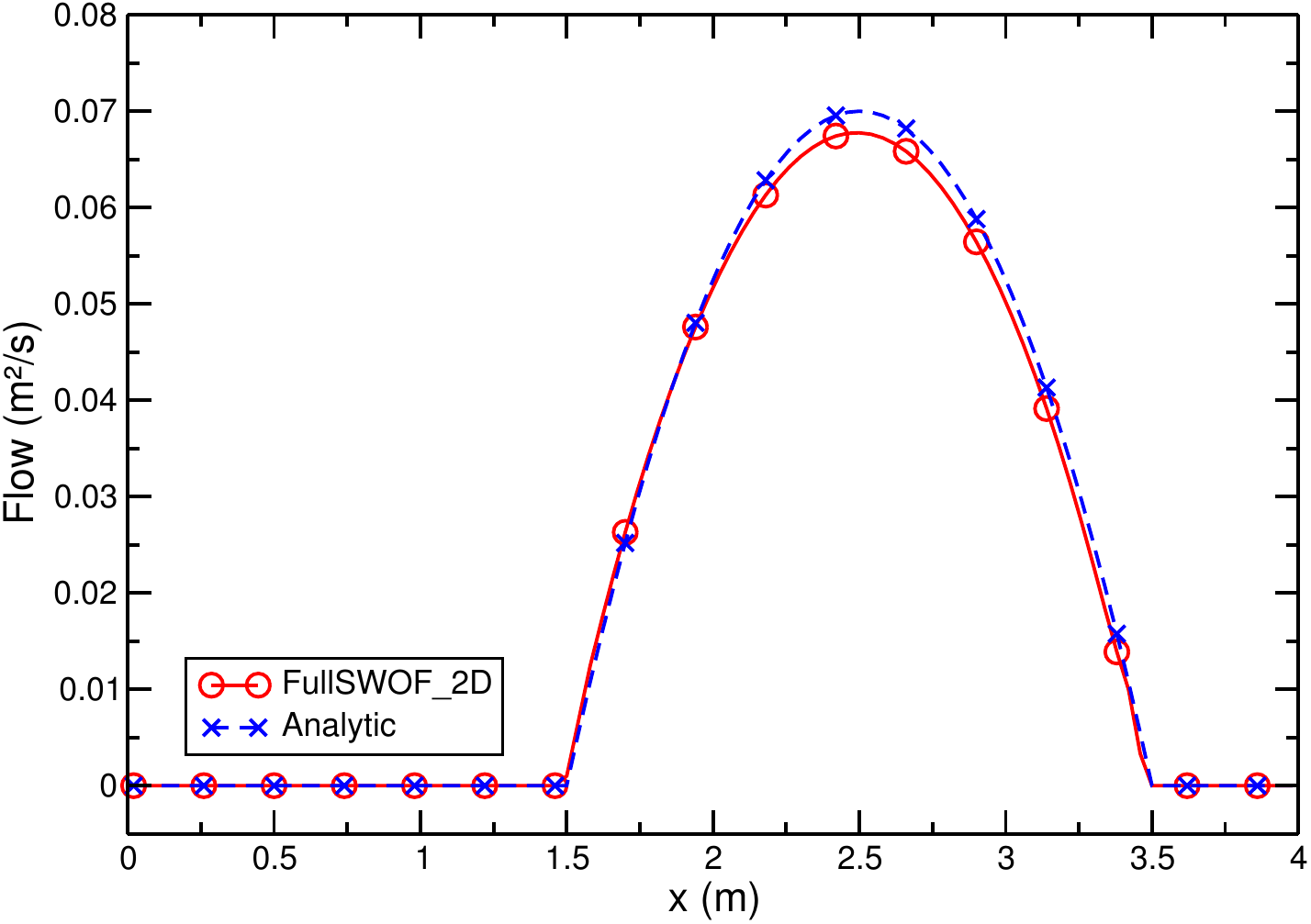}}
	\caption{Planar surface rotating in a paraboloid: 
		comparison of the analytic solution with the FullSWOF\_2D results.
		Cross-sections at $y~= 2\;\si{m}$.
		Simulation with 100~cells in~$x$ and~$y$ at T~= 13.4571~\si{s}.\label{Fig:Thacker}}
\end{figure}

Finally, we present a transitory two-dimensional solution:
a planar water surface rotating in a paraboloid \citep[§4.2.2]{SWASHES13}.
In the literature, this analytic solution is known as Thacker's 2D case, named after its author \citep{Thacker81}.
The topography is a paraboloid of revolution, and the shoreline is a moving circle.
The free surface exhibits periodic motion and remains planar in time.
To visualize this case, wine being swirled in a glass can be considered (cross-section in Figure~\ref{Fig:Thacker}a). 
This is a solution with a variable slope (in space) for which the wet-dry transitions are constantly moving. 
Because there is no friction, the rotation should not decrease over time.
Hence, this case tests both the ability of the schemes to simulate flows with comings and goings
and the numerical diffusion of the scheme (which causes a damping in the water depth over time).
The analytic solution at $t=0\;\si{s}$ is taken as the initial condition for the computation. The results are considered after
three periods.

Overall, FullSWOF\_2D produces results in good agreement with the analytic solution (Figure~\ref{Fig:Thacker}).
There is no spurious point at the wet-dry transitions, exemplifying
the capability of the well-balanced
schemes to properly compute these common situations for natural surface flows.
The slope of the water surface given by FullSWOF\_2D is slightly lower than that of the analytic solution.
In addition, the water flow rate is slightly lower than expected. This is a consequence of the numerical diffusion of the scheme.
Several techniques can be applied to improve the results of this specific case, but they may not be relevant for the general 
case because they would increase the complexity of the method and the computational cost.

\subsection{MacDonald-Type Solutions}\label{Sec:MacDo}
We turn now to a class of stationary solutions, more complex than the previous solutions (section~\ref{Sec:Classical}) in the sense that
the model is now complete, with friction and rain.
These solutions are
 obtained by the procedure introduced by I.~MacDonald: the water depth profile and the discharge are given, and the 
corresponding topography is then computed. In the original works,
the Manning friction law was considered \citep{MacDonald96, MacDonald97}, but the method allows many variants: other friction laws, rain, diffusion, etc.

The three selected solutions show the ability of the scheme to address
stationary states induced by the topography and by friction under a wide range
of flow conditions.

\subsubsection{Short Channel with a Smooth Transition and a Shock}

\begin{figure}[htbp]
	\subcaptionbox{Water depths and topography.\label{Fig:MacDoFluMan:h}}
		{\includegraphics[width=0.48\textwidth]{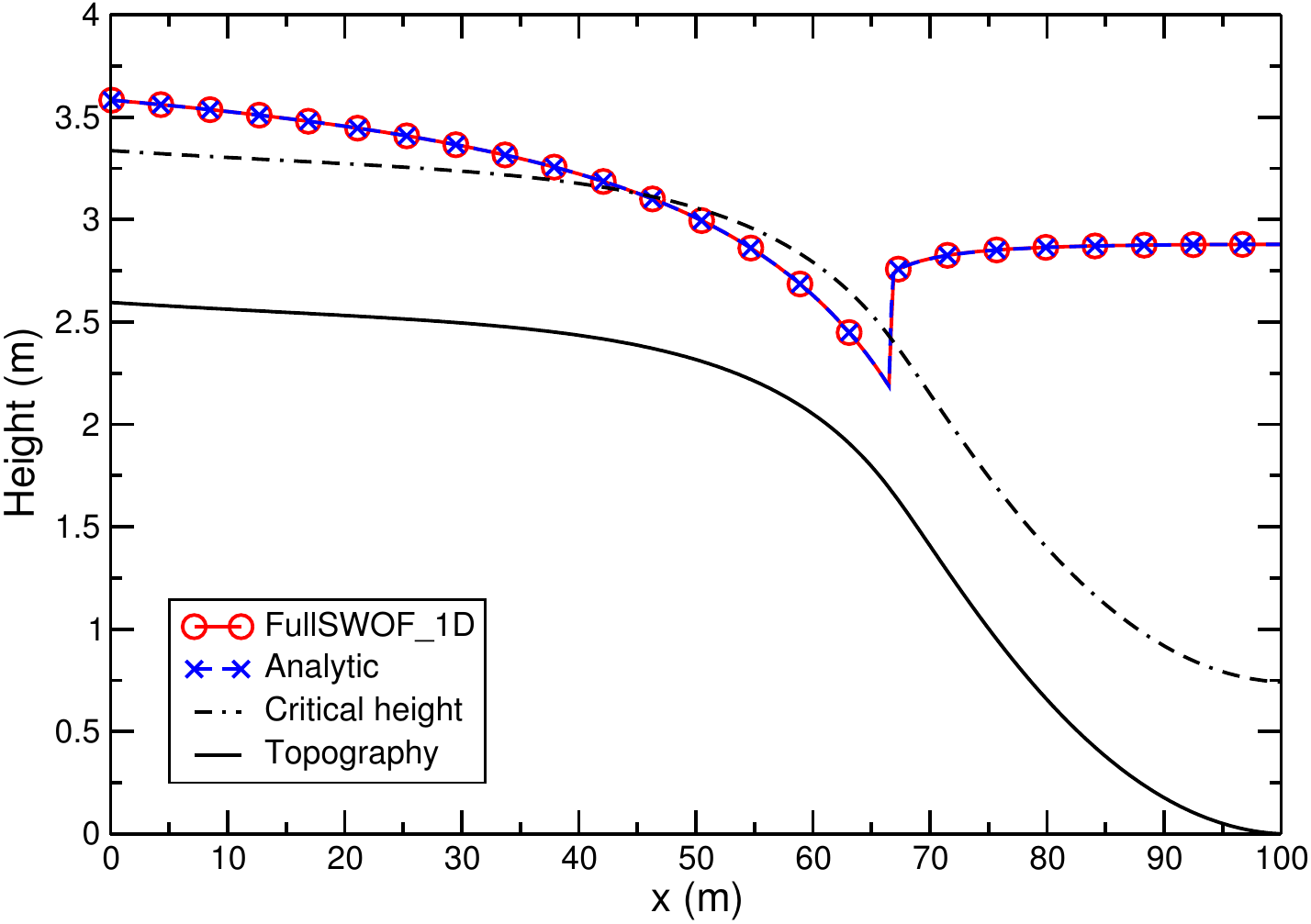}}
	\hfill
	\subcaptionbox{Water flow rates.\label{Fig:MacDoFluMan:q}}
		{\includegraphics[width=0.48\textwidth]{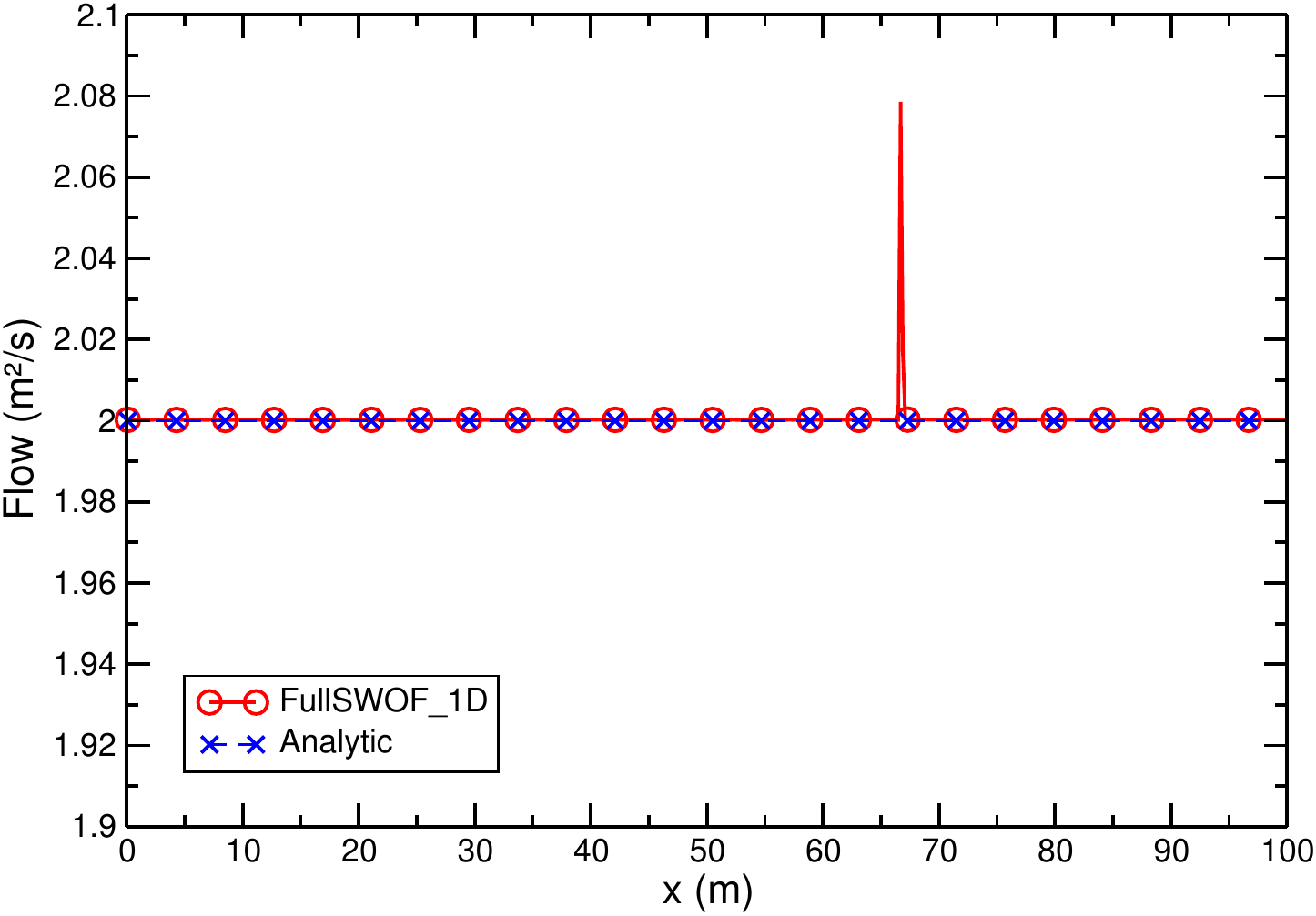}}
	\caption{Short channel with a smooth transition and a shock:
		comparison of the analytic solution with the FullSWOF\_1D results.
		Simulation with 500~cells at T~= 1500~\si{s}.\label{Fig:MacDoFluMan}}
\end{figure}

The length of the channel is $100$ m, and the discharge at steady state is $q=2\;\si{m^{2}.s^{-1}}$ (Figure~\ref{Fig:MacDoFluMan}). 
The flow is subcritical both upstream and downstream \citep[§3.2.2]{SWASHES13}.
In the intermediate part, the flow is supercritical. 
Hence, the flow goes from subcritical to supercritical via a sonic point and then ---
through a shock (located at $x=200/3 \approx 66.67\;\si{m}$) ---
becomes subcritical again (Figure~\ref{Fig:MacDoFluMan:h}). The Manning friction coefficient $n$ is equal to 
0.0328~\si{m^{-1/3}.s}. The run was performed at a resolution of 0.2~\si{m}.

The water flow rate simulated by FullSWOF\_1D matches the analytic solution except at the shock (Figure~\ref{Fig:MacDoFluMan:q}). A careful
examination shows that the mismatch is limited to two cells.
This mismatch is due to the product of approximations of two
discontinuous functions ($u$ and $h$) on a discretization grid.
Further testing showed that the shock was better represented with finer grids, which is the expected behavior. Future developments --- by us or
by others, considering that FullSWOF is open source --- could improve the results of this test case.

\subsubsection{Rain on a Long Channel with a Supercritical Flow}
\begin{figure}[htbp]
	\subcaptionbox{Water depths and topography.\label{Fig:Rain_tor:h}}
		{\includegraphics[width=0.48\textwidth]{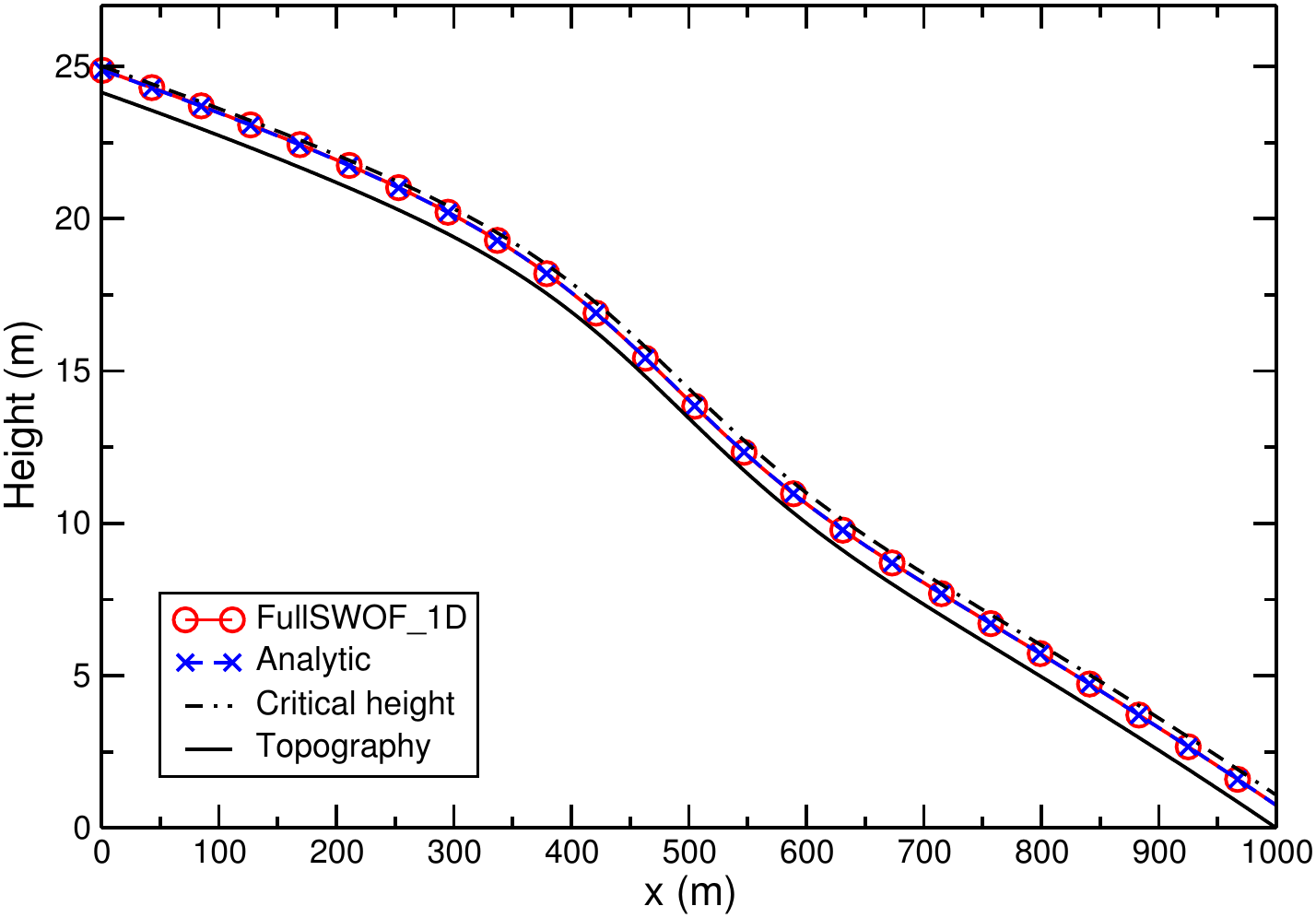}}
	\hfill
	\subcaptionbox{Water flow rates.\label{Fig:Rain_tor:q}}
		{\includegraphics[width=0.48\textwidth]{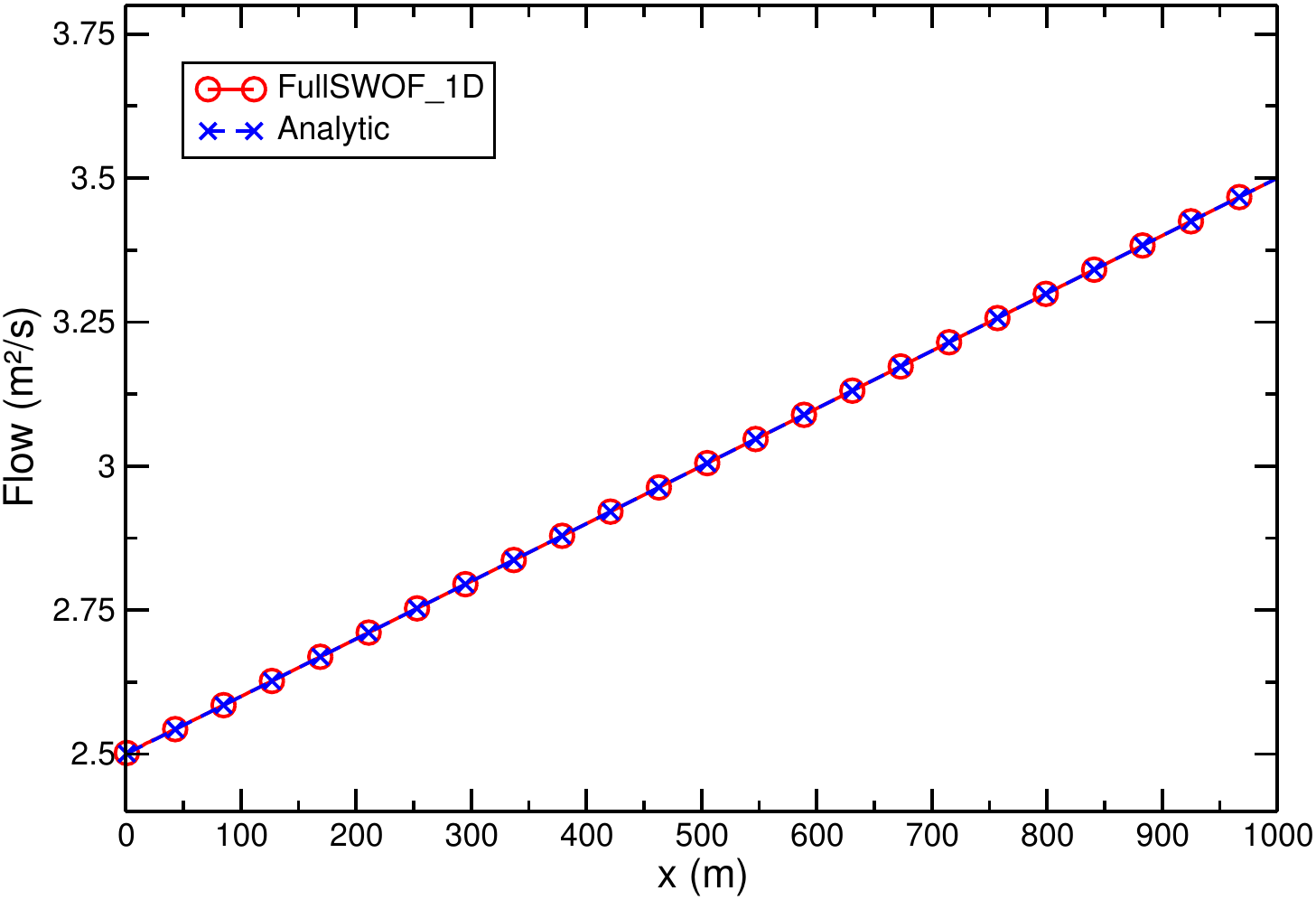}}
	\caption{Rain on a long channel with a supercritical flow:
		comparison of the analytic solution with the FullSWOF\_1D results.
		Simulation with 500~cells in~$x$ at T~= 3000~\si{s}.\label{Fig:Rain_tor}}
\end{figure}

This analytic case is similar to the previous case but includes rain
\citep[§3.3.2]{SWASHES13}. The channel length is set equal to 1000~\si{m} (Figure~\ref{Fig:Rain_tor}). 
Because the flow is supercritical along the entire channel, we consider a constant discharge and a constant water depth at inflow and a free outflow.
In the simulation, the channel is initially dry.
There is no rain until 1500~\si{s}; after this time, the rain
intensity is set to 0.001~\si{m.s^{-1}} until the end of the simulation.
The Darcy-Weisbach friction coefficient is $f=0.065$, and the inflow discharge is $q_0=2.5\;\si{m^{2}.s^{-1}}$ (Figure~\ref{Fig:Rain_tor}).

The results from FullSWOF\_1D are in good agreement with the analytic
solution, which demonstrates that the inclusion of the rain source term in the model and its implementation in the code are satisfactory.

\subsubsection{Pseudo-2D Channel with a Supercritical Flow}

In this section, we consider a pseudo-2D shallow-water system in a channel with a varying width.
By pseudo-2D, we mean an intermediate model between the one-dimensional and two-dimensional models.
More precisely, the pseudo-2D shallow-water equations are obtained by averaging all quantities both in the vertical direction 
and over the width of the channel ($y$ direction), i.e., perpendicular to the flow (recall that the two-dimensional shallow-water equations involve an average solely in the vertical direction).
The derivation is detailed in \citet{Goutal11}. In that paper, the authors showed that the new terms that appear in the pseudo-2D
shallow-water system due to the $y$ direction averaging procedure produce numerical difficulties that must be treated carefully. We perform a simulation
of the full two-dimensional equations with a topography and discharge corresponding to a stationary state for the pseudo-2D system, which 
enables some comparisons.

\begin{figure}[htbp]
	\subcaptionbox{Water depth simulated by FullSWOF\_2D. Inside
          the channel, the water depth ranges from 0.38 to 1.07~\si{m}.\label{Fig:MacDop2D:hdessus}}
		{\includegraphics[width=0.52\textwidth]{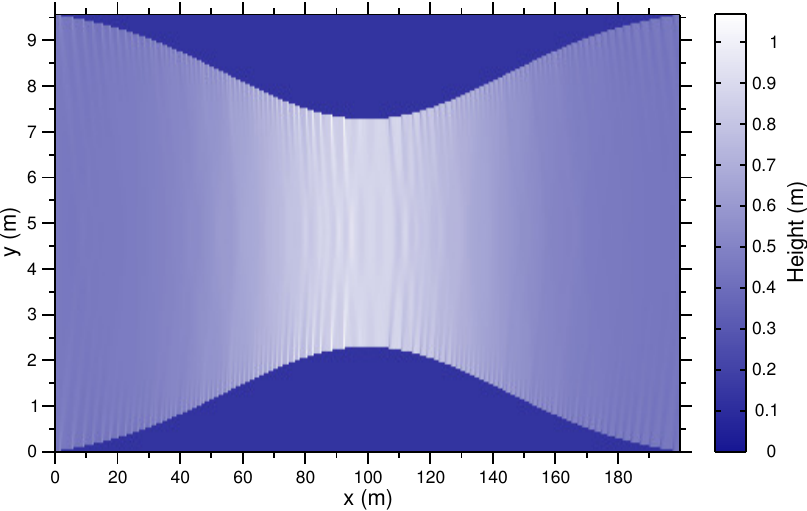}}
	\hfill
	\subcaptionbox{Water depths and topography. For FullSWOF\_2D, the mean depth inside the channel is plotted.\label{Fig:MacDop2D:h}}
		{\includegraphics[width=0.45\textwidth]{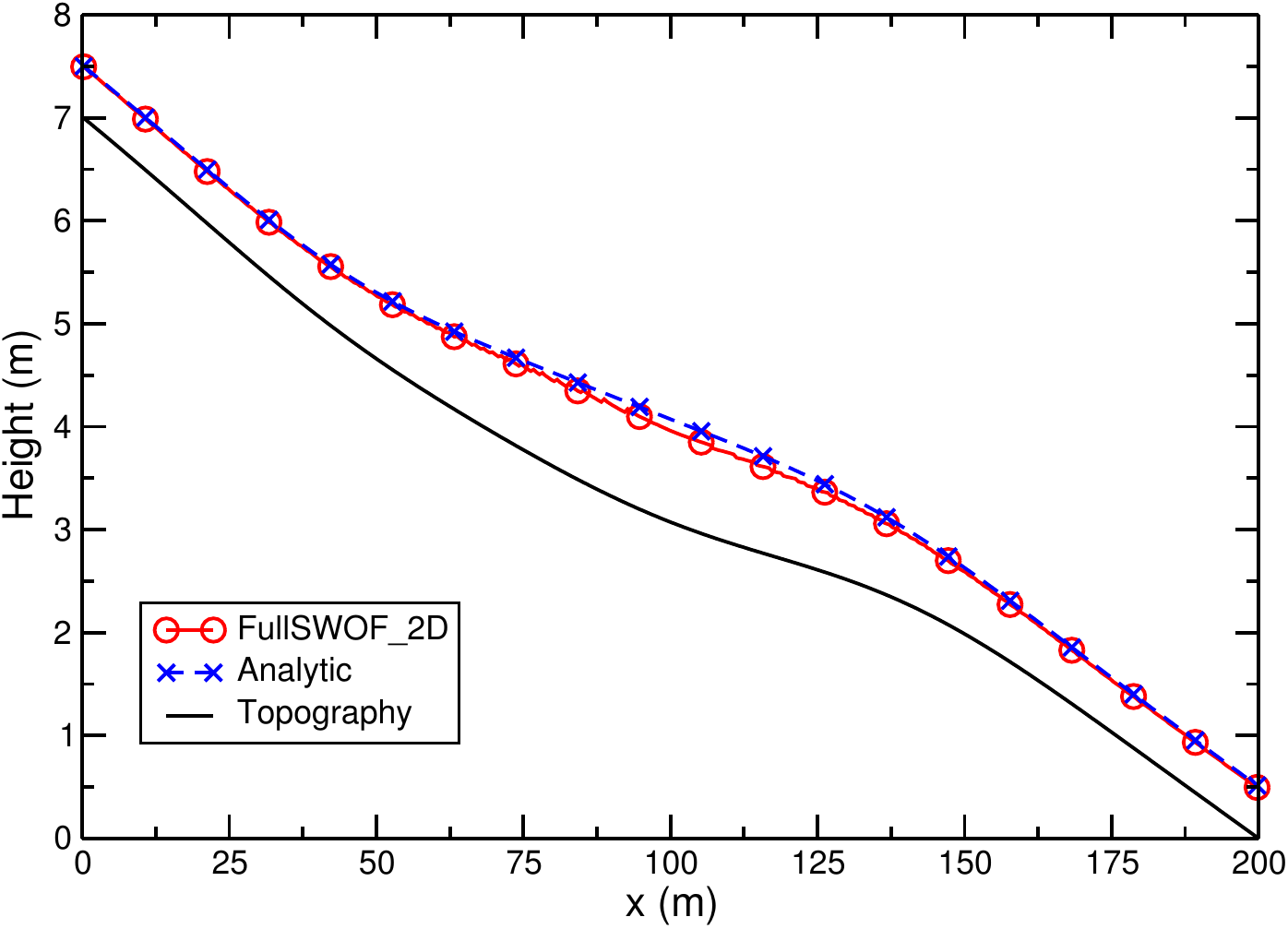}}
	\caption{Pseudo-2D channel with a supercritical flow:
		comparison of the analytic solution with the FullSWOF\_2D results.
		Simulation with 400~cells in~$x$ and 201~cells in~$y$ at T~= 200~\si{s}.\label{Fig:MacDop2D}}
\end{figure}

We consider a 200-\si{m}-long channel with a rectangular cross-section. The width and slope of the channel 
depend on $x$ (Figure~\ref{Fig:MacDop2D:hdessus}).
The inflow rate is fixed at $q=20\;\si{m^{3}.s^{-1}}$, and the water depth is prescribed at outflow \citep[§3.5.2]{SWASHES13}. The Manning coefficient is set to 0.03\;\si{m^{-1/3}.s}. The channel is initially dry, with a small puddle downstream (because of the outflow condition).

At steady state, FullSWOF\_2D produces a depth profile similar
to the analytic solution (Figure~\ref{Fig:MacDop2D:h}). An exact match
is not to be expected because the results of a 2D code are being compared with the analytic solution of a \emph{pseudo-}2D case. 

\subsection{Real Datasets}\label{Sec:Realdata}
Finally, we provide the results from FullSWOF using three real
datasets. The first case is an original laboratory experiment, the
second case comes from an experimental plot in Senegal, and the third
case is the well-known 
Malpasset dam break. We do not pretend to provide new insights into these
phenomena; we merely demonstrate the ability of FullSWOF 
to simulate real situations without data filtering. The refined
calibration of the parameters and a comparison with other pieces of software 
will be the topic of future works.

\subsubsection{Flow over a Corrugated Bottom}
\begin{figure}[htbp]
	\subcaptionbox{Close-up photography of the experimental setup. The red line is the measured profile.\label{Fig:Etape1:photo}}
		{\includegraphics[width=0.4\textwidth]{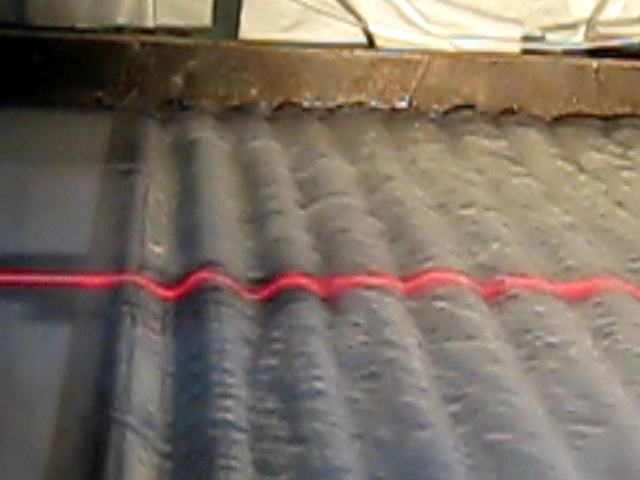}\\ \vspace{0.3em}}
	\hfill
	\subcaptionbox{Measured topography, experimental water surface and water surface computed by
		FullSWOF\_1D (Simulation with 1104~cells at T~= 10~\si{s}).\label{Fig:Etape1:profile}}
		{\includegraphics[width=0.5\textwidth]{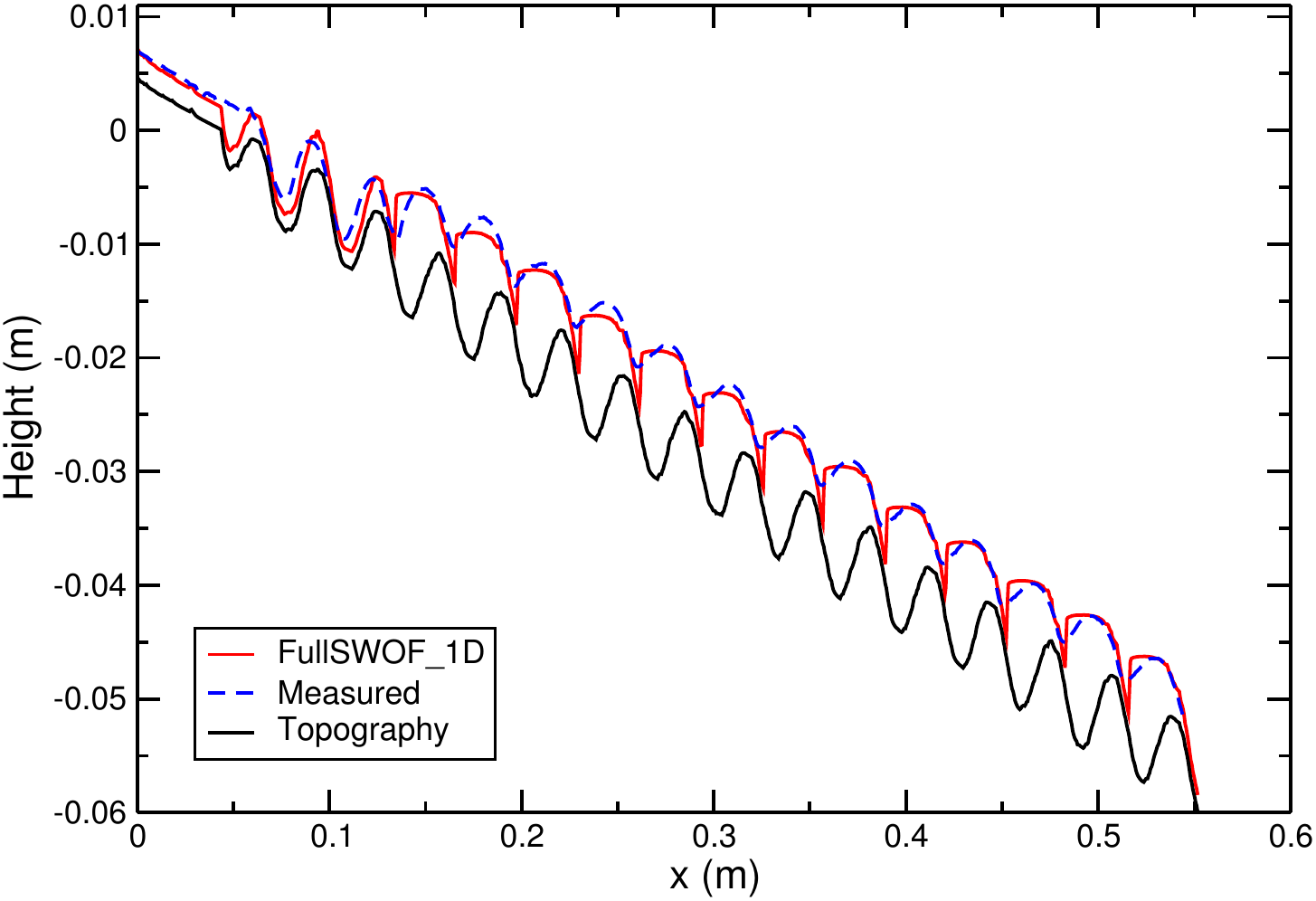}}
	\caption{Experimental flow over a corrugated bottom.\label{Fig:Etape1}}
\end{figure}

This first example is an experimental flow over an inclined channel having a corrugated bottom (Figure~\ref{Fig:Etape1}).
The discharge is imposed upstream (measured value:
0.69~\si{l.m^{-1}.s^{-1}}). The water
depths are measured at steady state along a 55-\si{cm}-long profile with a 0.5-\si{mm} resolution 
using the device described in \citet{Legout12}.
FullSWOF\_1D is run using the measured topography profile as the input. The only calibrated parameter is the Manning friction coefficient $n=0.0127\;\si{m^{-1/3}.s}$.

The comparison between the computed and measured values shows that 
FullSWOF\_1D is able to reproduce the qualitative behavior
after the first bump ($x > 0.06\;\si{m}$)  (Figure~\ref{Fig:Etape1:profile}). 
FullSWOF\_1D reproduces all the hydraulic jumps (for $x >
0.13\;\si{m}$) and correctly locates them at the minima of the measured solution.
However, the simulated solution exhibits shocks that are steeper than the measured shocks.
This, as well as the poor restitution at the beginning of the channel
(at approximately $x = 0.05\;\si{m}$), is likely due to a shortcoming of the model.
Indeed, in this range of water depths and velocities, the surface tension likely cannot be neglected. This is a possible 
extension of the FullSWOF software.

\subsubsection{Rain over a Field Plot (Thi\`es, Senegal)}\label{Sec:Thies}

The purpose of this section is to test FullSWOF\_2D using a real system: a plot in
Thi\`es, Senegal \citep{Tatard08, Mugler11}, an experimental system instrumented by IRD\@.
This system consists of an artificial rainfall of approximately 2
hours with an intensity of approximately 70~\si{mm.h^{-1}}  on a
sandy-soil plot with an area of $4 \times 10\;\si{m^2}$. The
plot has the classical configuration of Wooding's open book, with a
1\% slope along the $Ox-$ and $Oy-$axes. 
The complete dataset is freely available at \url{http://www.umr-lisah.fr/Thies_2004/}. 
FullSWOF\_2D is used with a variable time step (the CFL value is fixed at 0.4) and with the same parameters as in \citet{Tatard08}: $f=0.26$, $h_f =0.06\;\si{m}$, $\theta_s-\theta_i =0.12$, $K_s=4.4\ 10^{-6}\;\si{m.s^{-1}}$ and $K_c=0$. 

\begin{figure}[htbp]
	\subcaptionbox{Topography (without the 1\% slope along~$y$).\label{Fig:ThiesFine:Z}}
		{\includegraphics[height=0.4\textheight]{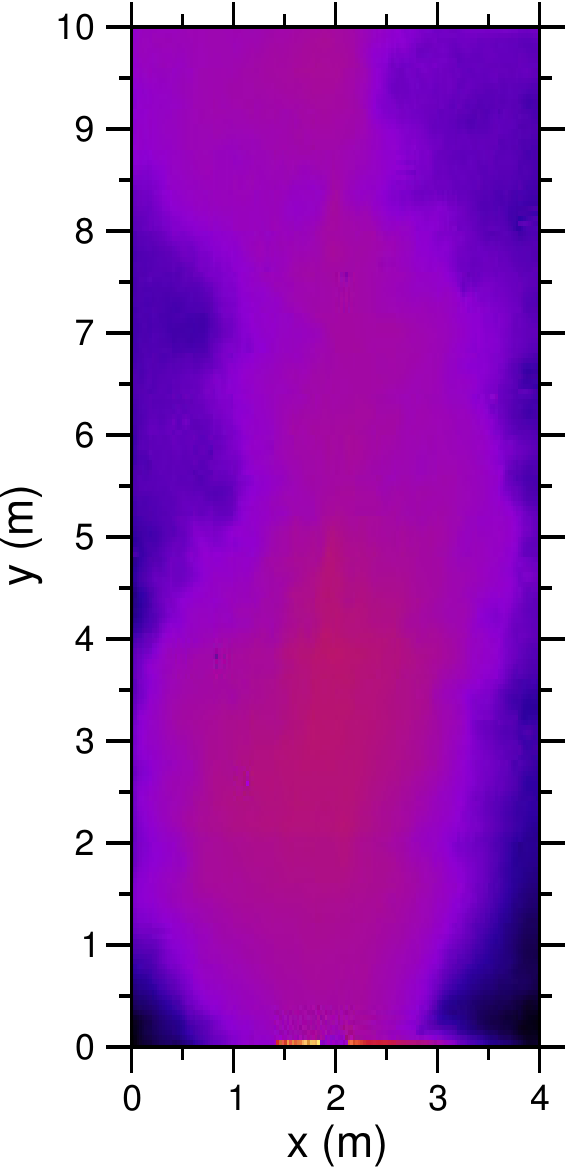}}
	\raisebox{3em}{\hspace*{1em}\includegraphics[width=0.15\textwidth]{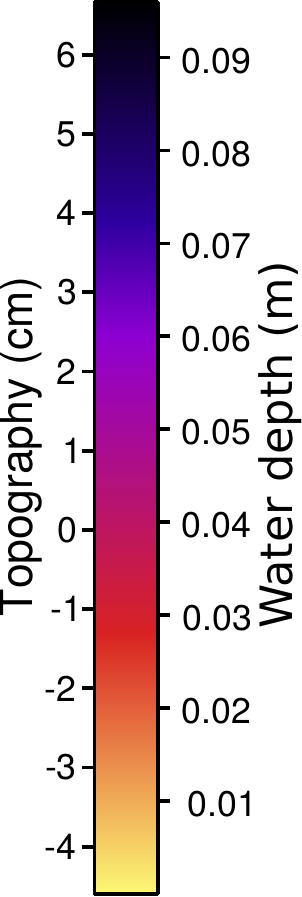}\hspace*{1em}}
	\subcaptionbox{Water depth.\label{Fig:ThiesFine:H}}
		{\includegraphics[height=0.4\textheight]{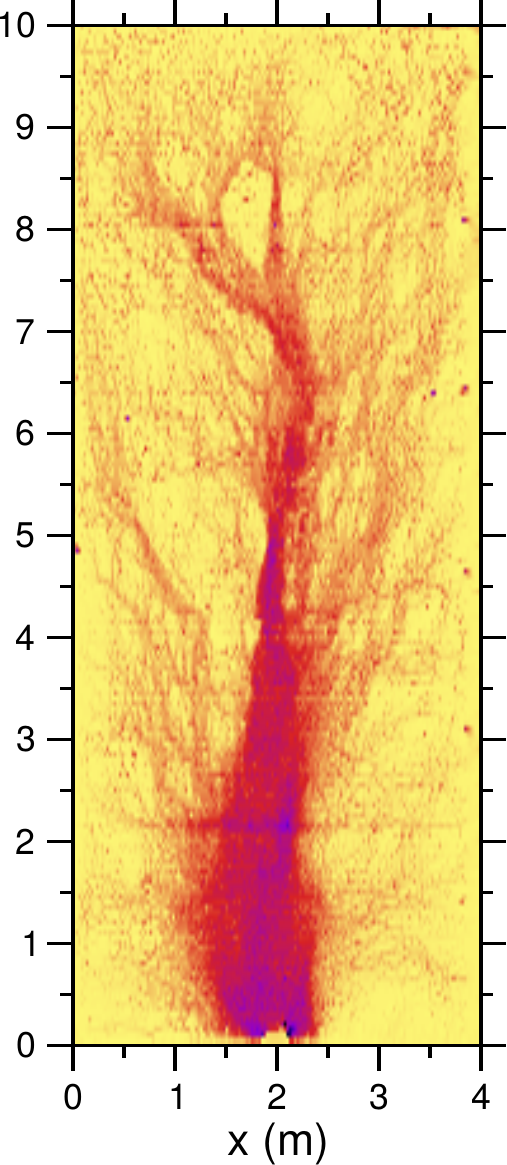}}
		
	\bigskip
		
	\subcaptionbox{Velocity.\label{Fig:ThiesFine:U}}
		{\includegraphics[height=0.4\textheight]{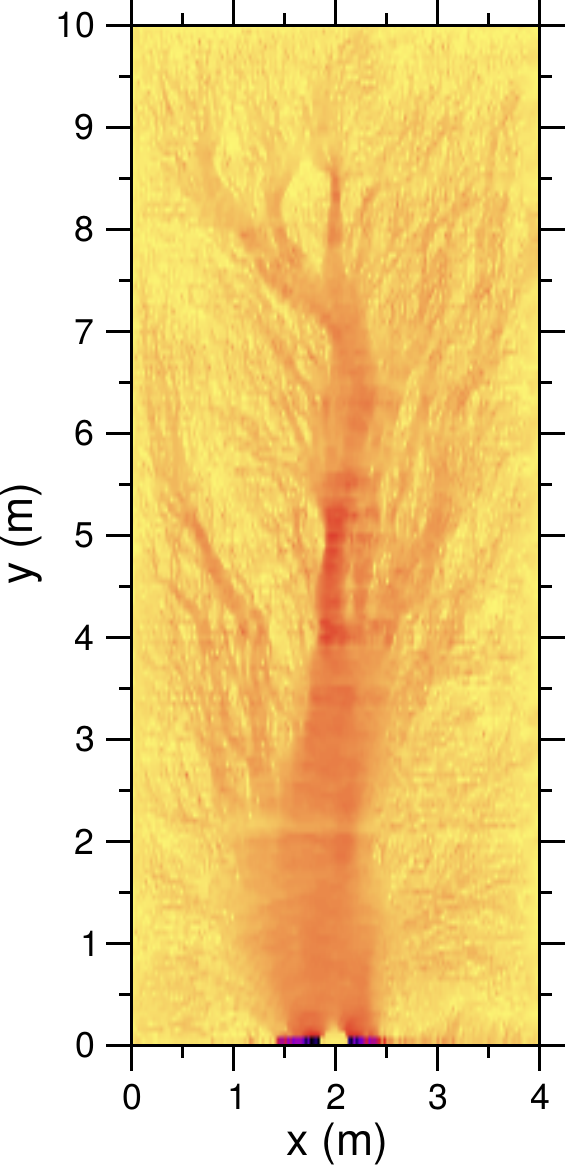}}
	\raisebox{3em}{\hspace*{1em}\includegraphics[width=0.15\textwidth]{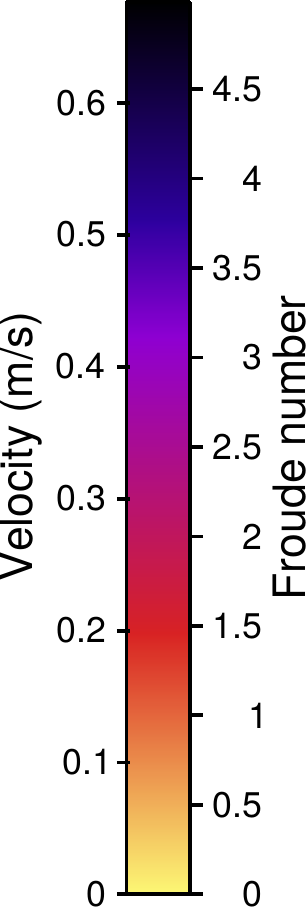}\hspace*{1em}}
	\subcaptionbox{Froude number.\label{Fig:ThiesFine:F}}
		{\includegraphics[height=0.4\textheight]{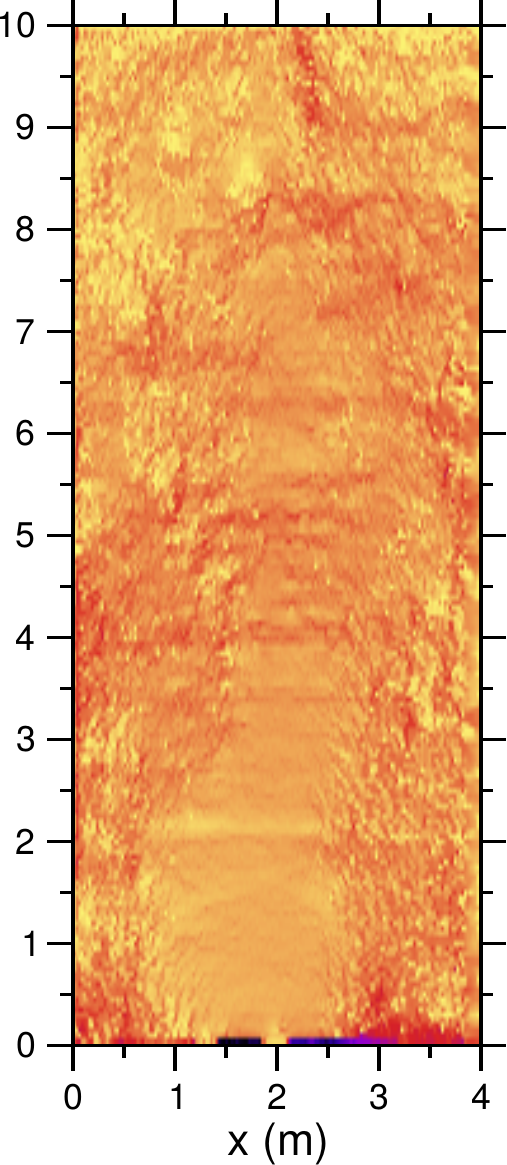}}
	\caption{Results of FullSWOF\_2D for the Thi\`es field plot.\\
		Simulation with 160~cells in~$x$ and 200~cells in~$y$ at T~= 7400~\si{s}.\label{Fig:ThiesFine}}
\end{figure}

Qualitatively, the results (Figure~\ref{Fig:ThiesFine}) are similar to
those obtained with other software packages \citep{Tatard08, Mugler11}.
At this stage, these results illustrate the ability of FullSWOF\_2D to simulate a dynamic flow without any
filtering of the data, in contrast to the previously utilized software
packages. 
A comparison will be performed in a more detailed study.

\subsubsection{The Malpasset Dam Break: a Large Dataset Demanding the Parallel Version of FullSWOF\_2D}

This section is devoted to a large-scale test case.
We implement the well-known Malpasset dam break, which occurred in 1959 in the south of France (\citet{Hervouet00},
 \citet[p.281-288]{Hervouet07}). Because of its varying topography and complex geometry, this is a classical test for numerical methods and
 hydraulics software validation.
For this case, we used the parallel version of FullSWOF\_2D
\citep{CEMRACS13} to reduce the computational time to a few hours.
The dimensions of the computing domain are $L_x=17273.9\;\si{m}$ with 1000~cells
and $L_y=8381.3\;\si{m}$ with 486~cells.
The topography comes from a map published before the disaster that
was subsequently digitized.
We consider the Manning law with $n=0.033\;\si{m^{-1/3}.s}$, as
advised in \citet{Hervouet00}. No calibration was performed for this
simulation, which uses a total simulation time of T~= 4000~\si{s}.

\begin{figure}[htbp]
	\subcaptionbox{t =   0 s\label{Fig:Malpasset:0000}}
		{\includegraphics[width=0.22\textwidth]{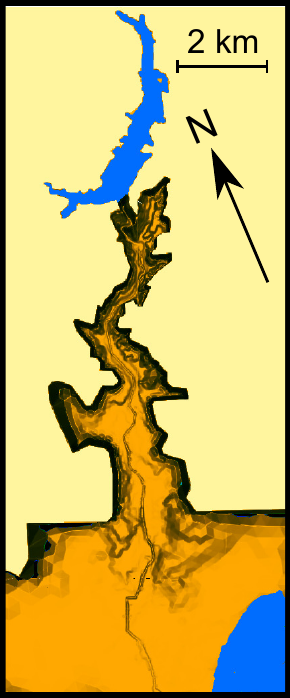}}
	\subcaptionbox{t =1000 s.\label{Fig:Malpasset:1000}}
		{\includegraphics[width=0.22\textwidth]{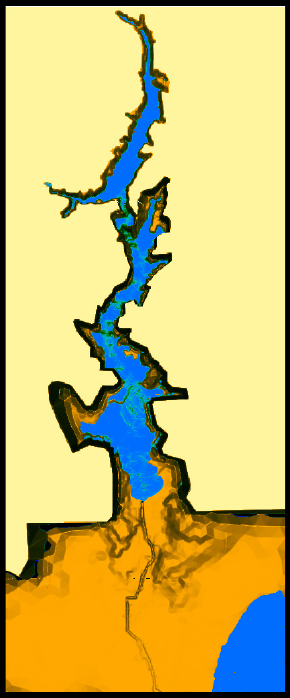}}
	\subcaptionbox{t =2500 s.\label{Fig:Malpasset:2500}}
		{\includegraphics[width=0.22\textwidth]{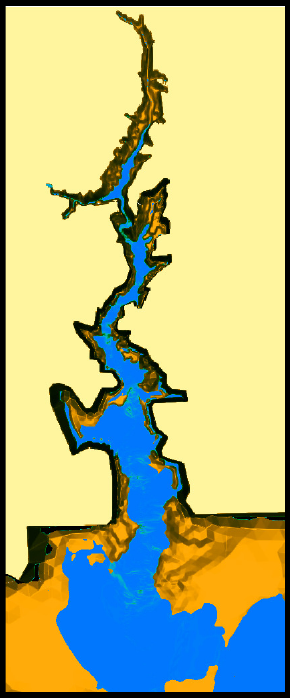}}
	\subcaptionbox{t = 4000 s.\label{Fig:Malpasset:4000}}
		{\includegraphics[width=0.22\textwidth]{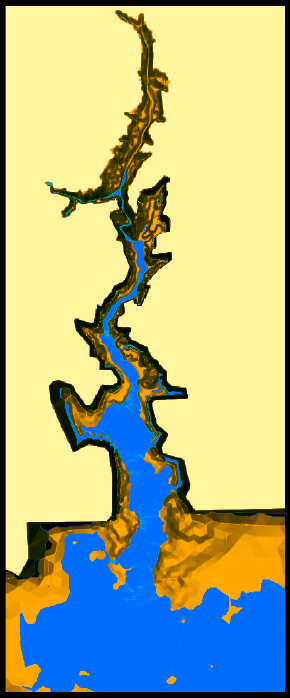}}\\
		\subcaptionbox{Gauge location.\label{Fig:Malpasset:Location}}
		{\includegraphics[width=0.22\textwidth]{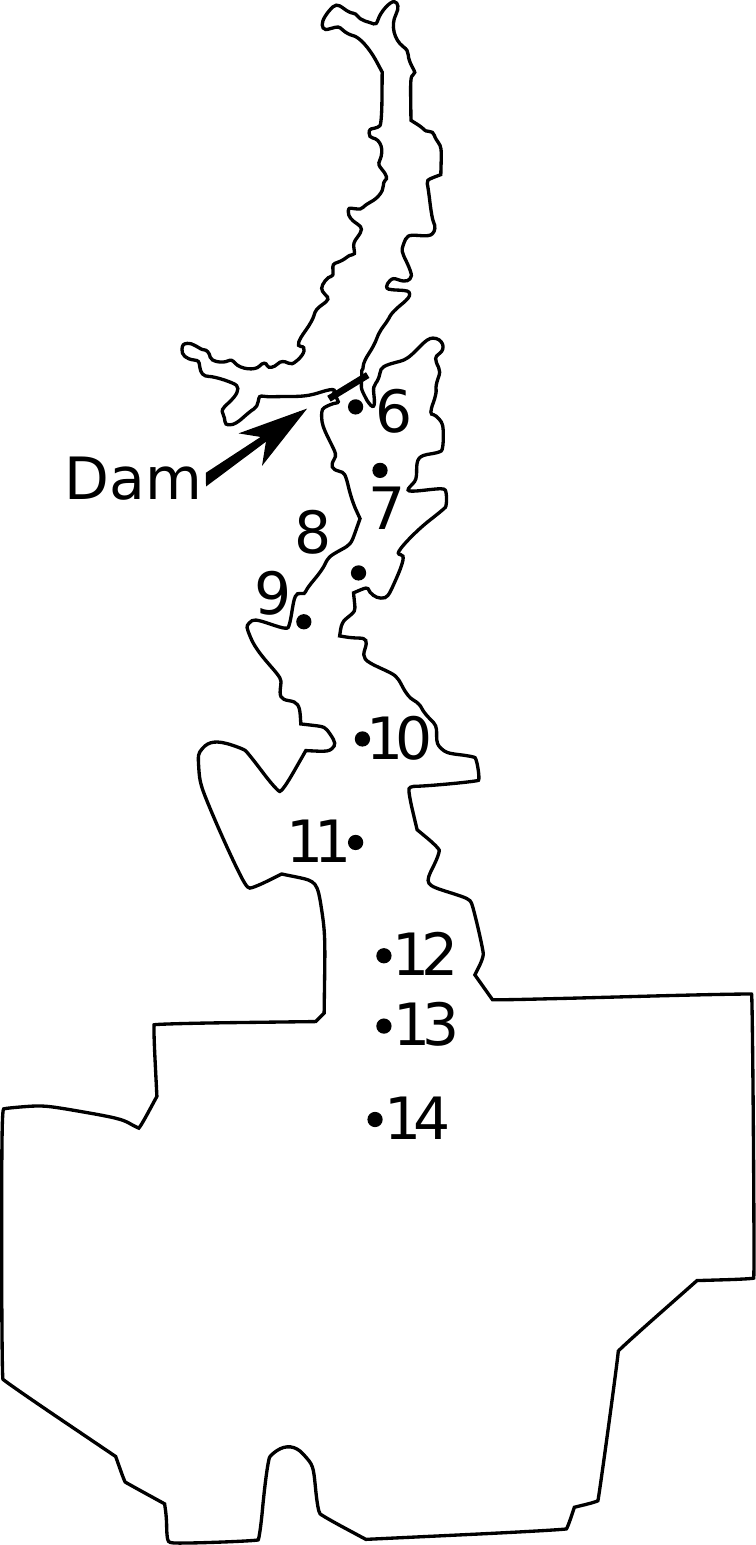}}\hspace{1em}
	\subcaptionbox{Comparison of maximal water elevations\label{Fig:Malpasset:Compare}}
		{\includegraphics[width=0.66\textwidth]{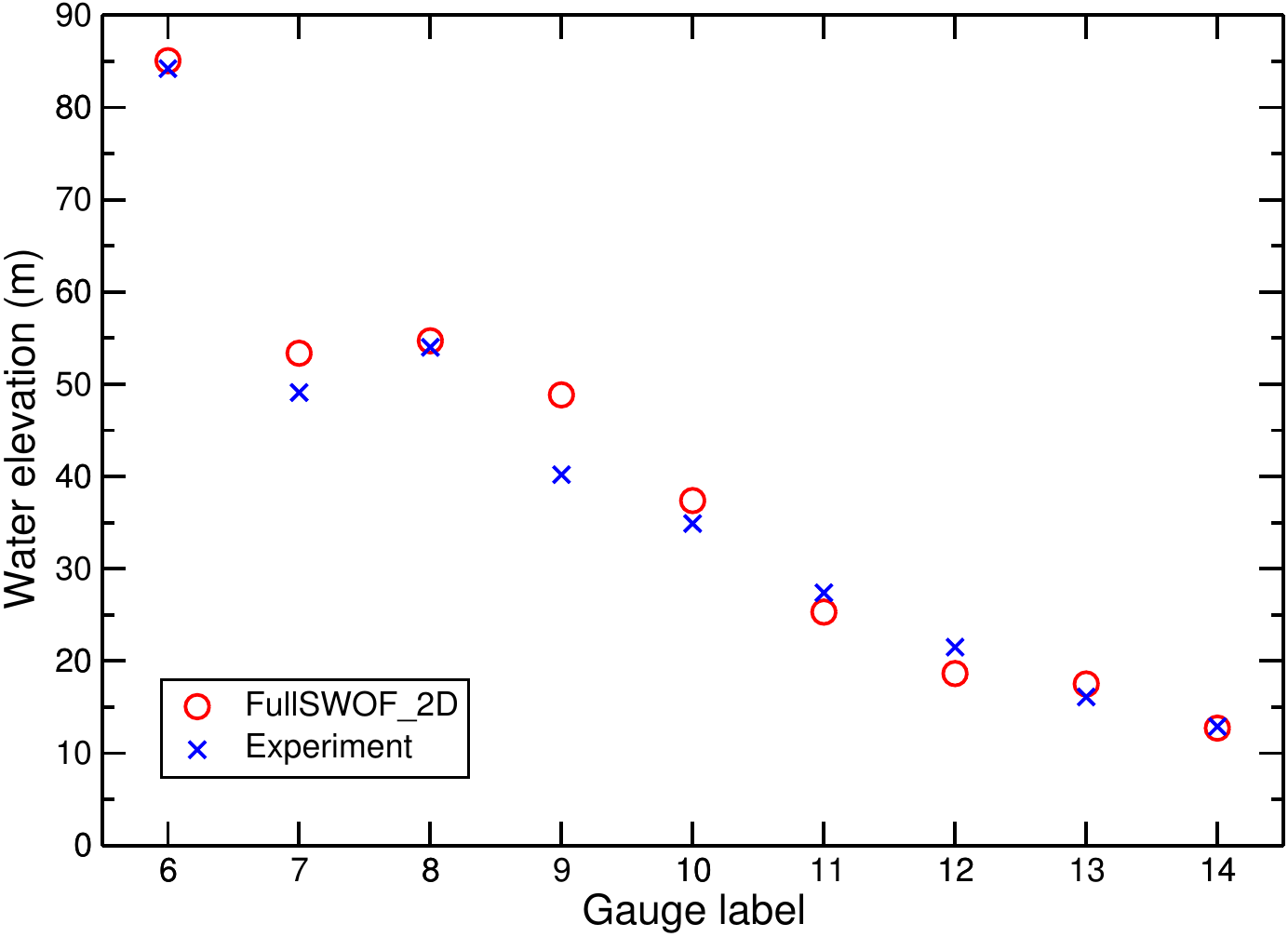}}
	\caption{Malpasset dam break: evolution of the water flow simulated by FullSWOF\_2D (parallel version) and comparison with experimental data.\label{Fig:Malpasset}}
\end{figure}

The results are presented in Figure~\ref{Fig:Malpasset}~a--d at four time steps.
A scaled physical experiment was built by the Laboratoire National d'Hydraulique in
1964 to study the dam-break flow (for more details, see \citet{Hervouet00} and \citet{Hervouet07}).
The maximum water level was recorded at nine gauges (labeled 6 to 14) during the physical experiment
(Figure~\ref{Fig:Malpasset}e).
The maximum water elevations of the scaled experiment are  well captured by FullSWOF\_2D (Figure~\ref{Fig:Malpasset}f).

\section{Conclusions}

The FullSWOF software project grew out of a long-term collaboration between mathematicians and hydrologists in Orléans (France). Several specific features make it particularly suitable for
applications in hydrology: FullSWOF considers rain, infiltration and
classical friction laws and can operate with DEMs. Moreover, the
software package has
demonstrated its ability to address a wide range of flow
conditions. The latest releases of the three codes (FullSWOF\_1D,
FullSWOF\_2D and FullSWOF\_UI) are fairly stable and robust. They are
freely available for users outside of the original team for industrial, scientific and educational purposes.

The project is still undergoing several improvements. Apart from the
usual bug correction work, several improvements to the current codes
(e.g., implementing non-homogeneous friction) are regularly produced, and a
parallelized 2D version is under development. More conceptual
modifications are under consideration, such as the introduction of
erosion and sedimentation models. Another direction of improvement
concerns the surface tension. Although these modifications currently  originate mainly from the initial development team, we emphasize that external contributions are welcome.

\section*{Acknowledgments}

This study is part of the ANR METHODE granted by the French National Agency for Research, ANR-07-BLAN-0232.
The authors would like to thank Olivier Planchon (IRD) for the data used in section~\ref{Sec:Thies} for the rain over a field plot (Thi\`es, Senegal).

\bibliographystyle{apalike}

\bibliography{FullSWOF}

\begin{thebibliography}{}

\bibitem[Audusse et~al., 2004]{Audusse04c}
Audusse, E., Bouchut, F., Bristeau, M.-O., Klein, R., and Perthame, B. (2004).
\newblock A fast and stable well-balanced scheme with hydrostatic
  reconstruction for shallow water flows.
\newblock {\em SIAM J. Sci. Comput.}, 25(6):2050--2065.

\bibitem[Baartman et~al., 2013]{Baartman13}
Baartman, J. E.~M., Temme, A. J. A.~M., Veldkamp, T., Jetten, V.~G., and
  Schoorl, J.~M. (2013).
\newblock Exploring the role of rainfall variability and extreme events in
  long-term landscape development.
\newblock {\em Catena}, 109:25--38.

\bibitem[Barnes, 2010]{Barnes10}
Barnes, N. (2010).
\newblock Publish your computer code: it is good enough.
\newblock {\em Nature}, 467(7317):753.

\bibitem[Barr{\'e}~de Saint-Venant, 1871]{saintvenant71}
Barr{\'e}~de Saint-Venant, A.-J.-C. (1871).
\newblock Th\'eorie du mouvement non-permanent des eaux, avec application aux
  crues des rivi\`eres et \`a l'introduction des mar\'ees dans leur lit.
\newblock {\em Comptes Rendus de l'Acad\'emie des Sciences}, 73:147--154.

\bibitem[Batten et~al., 1997]{Batten97}
Batten, P., Clarke, N., Lambert, C., and Causon, D.~M. (1997).
\newblock On the choice of wavespeeds for the {HLLC} {R}iemann solver.
\newblock {\em SIAM J. Sci. Comput.}, 18(6):1553--1570.

\bibitem[Berger et~al., 2011]{GeoClaw11}
Berger, M.~J., George, D.~L., LeVeque, R.~J., and Mandli, K.~T. (2011).
\newblock The {G}eo{C}law software for depth-averaged flows with adaptive
  refinement.
\newblock {\em Advances in Water Resources}, 34:1195--1206.

\bibitem[Berm\'udez and V\'azquez, 1994]{Bermudez94}
Berm\'udez, A. and V\'azquez, M.~E. (1994).
\newblock Upwind methods for hyperbolic conservation laws with source terms.
\newblock {\em Computers \& Fluids}, 23(8):1049--1071.

\bibitem[Bouchut, 2004]{Bouchut04}
Bouchut, F. (2004).
\newblock {\em Nonlinear Stability of Finite Volume Methods for Hyperbolic
  Conservation Laws, and Well-Balanced Schemes for Sources}.
\newblock Frontiers in Mathematics. Birkh{\"a}user Basel.

\bibitem[Bristeau and Coussin, 2001]{Bristeau01}
Bristeau, M.-O. and Coussin, B. (2001).
\newblock Boundary conditions for the shallow water equations solved by kinetic
  schemes.
\newblock Technical Report RR-4282, INRIA.

\bibitem[Brunner, 2010]{HECRAS-Hyd11}
Brunner, G.~W. (2010).
\newblock {\em HEC-RAS River Analysis System. Hydraulic Reference Manual
  Version 4.1}.
\newblock US Army Corps of Engineers. Hydrologic Engineering Center, Davis, CA,
  USA, cpd-69 edition.

\bibitem[Chow, 1959]{Chow59}
Chow, V.~T. (1959).
\newblock {\em Open-Channel Hydraulics}.
\newblock Mc{G}raw-Hill, New York.

\bibitem[Claerbout and Karrenbach, 1992]{Claerbout92}
Claerbout, J. and Karrenbach, M. (1992).
\newblock Electronic documents give reproducible research a new meaning.
\newblock In {\em Proceedings of the 62$^{nd}$ Annual International Meeting of
  the Society of Exploration Geophysics}, pages 601--604.

\bibitem[Cordier et~al., 2013]{CEMRACS13}
Cordier, S., Coullon, H., Delestre, O., Laguerre, C., Le, M.~H., Pierre, D.,
  and Sadaka, G. (2013).
\newblock {FullSWOF\_Paral}: Comparison of two parallelization strategies
  ({MPI} and {SkelGIS}) on a software designed for hydrology applications.
\newblock {\em ESAIM Proceeding. CEMRACS~2012: Numerical Methods and Algorithms
  for High Performance Computing}, 44:59--79.

\bibitem[Delestre, 2010]{Delestre10b}
Delestre, O. (2010).
\newblock {\em Rain water overland flow on agricultural fields simulation}.
\newblock PhD thesis, University of Orl\'eans, France.
\newblock In French.

\bibitem[Delestre et~al., 2014]{Delestre14a}
Delestre, O., Cordier, S., Darboux, F., Du, M., James, F., Laguerre, C., Lucas,
  C., and Planchon, O. (2014).
\newblock {FullSWOF}: a software for overland flow simulation.
\newblock In Gourbesville, P., Cunge, J., and Caignaert, G., editors, {\em
  Advances in Hydroinformatics. SIMHYDRO 2012 – New Frontiers of Simulation},
  pages 221--231, Nice, France. Société Hydrotechnique de France (SHF), the
  University of Nice-Sophia Antipolis (UNS), the International Association for
  Hydraulic Research (AIHR) and the Association Française de Mécanique (AFM),
  Springer.

\bibitem[Delestre et~al., 2012]{CRAShydro}
Delestre, O., Cordier, S., Darboux, F., and James, F. (2012).
\newblock A limitation of the hydrostatic reconstruction technique for
  {S}hallow {W}ater equations.
\newblock {\em Comptes Rendus de l'Acad\'emie des Sciences - S\'erie I -
  Math\'ematique}, 350:677--681.

\bibitem[Delestre et~al., 2009]{Delestre09}
Delestre, O., Cordier, S., James, F., and Darboux, F. (2009).
\newblock Simulation of rain-water overland-flow.
\newblock In Tadmor, E., Liu, J.-G., and Tzavaras, A., editors, {\em
  Proceedings of the $\rm 12^{th}$ International Conference on Hyperbolic
  Problems}, volume~67 of {\em Proceedings of Symposia in Applied Mathematics},
  pages 537--546, University of Maryland, College Park (USA). Amer. Math. Soc.

\bibitem[Delestre and James, 2010]{DelestreJames09}
Delestre, O. and James, F. (2010).
\newblock Simulation of rainfall events and overland flow.
\newblock In {\em Proceedings of 10$^{th}$ International Conference
  {Z}aragoza-{P}au on Applied Mathematics and Statistics, {J}aca, {S}pain},
  volume~35 of {\em Monograf\'ias Matem\'aticas Garc\'ia de Galdeano}, pages
  125--135.

\bibitem[Delestre et~al., 2013]{SWASHES13}
Delestre, O., Lucas, C., Ksinant, P.-A., Darboux, F., Laguerre, C., Vo, T.
  N.~T., James, F., and Cordier, S. (2013).
\newblock {SWASHES}: a compilation of {S}hallow-{W}ater analytic solutions for
  hydraulic and environmental studies.
\newblock {\em International Journal for Numerical Methods in Fluids},
  72:269--300.

\bibitem[{DHI Software}, 2009]{MIKE11}
{DHI Software} (2009).
\newblock {\em {MIKE}~11. A Modelling System for Rivers and Channels. Reference
  Manual}.
\newblock Denmark.

\bibitem[Esteves et~al., 2000]{Esteves00}
Esteves, M., Faucher, X., Galle, S., and Vauclin, M. (2000).
\newblock Overland flow and infiltration modelling for small plots during
  unsteady rain: numerical results versus observed values.
\newblock {\em Journal of Hydrology}, 228(3--4):265--282.

\bibitem[Fiedler and Ramirez, 2000]{Fiedler00}
Fiedler, F.~R. and Ramirez, J.~A. (2000).
\newblock A numerical method for simulating discontinuous shallow flow over an
  infiltrating surface.
\newblock {\em International Journal for Numerical Methods in Fluids},
  32(2):219--239.

\bibitem[Godlewski and Raviart, 1996]{Godlewski96}
Godlewski, E. and Raviart, P.-A. (1996).
\newblock {\em Numerical approximation of hyperbolic systems of conservation
  laws}, volume 118 of {\em Applied Mathematical Sciences}.
\newblock Springer-Verlag, New York, USA.

\bibitem[Goutal and Maurel, 1997]{Goutal97}
Goutal, N. and Maurel, F. (1997).
\newblock Proceedings of the $\rm 2^{nd}$ workshop on dam-break wave
  simulation.
\newblock Technical Report HE-43/97/016/B, Electricit\'{e} de France, Direction
  des \'{e}tudes et recherches.

\bibitem[Goutal and Sainte-Marie, 2011]{Goutal11}
Goutal, N. and Sainte-Marie, J. (2011).
\newblock A kinetic interpretation of the section-averaged {S}aint-{V}enant
  system for natural river hydraulics.
\newblock {\em International Journal for Numerical Methods in Fluids},
  67(7):914--938.

\bibitem[Green and Ampt, 1911]{GreenAmpt11}
Green, W.~H. and Ampt, G.~A. (1911).
\newblock Studies on soil physics. part~{I}.---the flow of air and water
  through soils.
\newblock {\em The Journal of Agricultural Science}, 4(1):1--24.

\bibitem[Greenberg and Le{R}oux, 1996]{Greenberg96}
Greenberg, J.~M. and Le{R}oux, A.~Y. (1996).
\newblock A well-balanced scheme for the numerical processing of source terms
  in hyperbolic equation.
\newblock {\em SIAM Journal on Numerical Analysis}, 33:1--16.

\bibitem[Halcrow, 2012]{ISIS}
Halcrow (2012).
\newblock {\em {ISIS} {2D} quick start guide --- v.~3.6}.

\bibitem[Harten et~al., 1983]{Harten83}
Harten, A., Lax, P.~D., and van Leer, B. (1983).
\newblock On upstream differencing and {G}odunov-type schemes for hyperbolic
  conservation laws.
\newblock {\em SIAM Review}, 25(1):35--61.

\bibitem[Harvey and Han, 2002]{Harvey02}
Harvey, H. and Han, D. (2002).
\newblock The relevance of open source to hydroinformatics.
\newblock {\em Journal of Hydroinformatics}, 4:219--234.

\bibitem[Hervouet, 2000]{Hervouet00}
Hervouet, J.-M. (2000).
\newblock A high resolution 2-{D} dam-break model using parallelization.
\newblock {\em Hydrological Processes}, 14(13):2211--2230.

\bibitem[Hervouet, 2007]{Hervouet07}
Hervouet, J.-M. (2007).
\newblock {\em Hydrodynamics of free surface flows}.
\newblock Wiley.

\bibitem[Hillel and Gardner, 1970]{Hillel70}
Hillel, D. and Gardner, W.~R. (1970).
\newblock Transient infiltration into crust-topped profiles.
\newblock {\em Soil Science}, 109(2):69--76.

\bibitem[Kolditz et~al., 2008]{Kolditz08}
Kolditz, O., Delfs, J.-O., Bürger, C., Beinhorn, M., and Park, C.-H. (2008).
\newblock Numerical analysis of coupled hydrosystems based on an
  object-oriented compartment approach.
\newblock {\em Journal of Hydroinformatics}, 10:227--244.

\bibitem[Kutija and Murray, 2007]{Kutija07}
Kutija, V. and Murray, M.~G. (2007).
\newblock An object-oriented approach to the modelling of free-surface flows.
\newblock {\em Journal of Hydroinformatics}, 9:81--94.

\bibitem[Lee and Wright, 2010]{Lee10}
Lee, S.-H. and Wright, N.~G. (2010).
\newblock Simple and efficient solution of the shallow water equations with
  source terms.
\newblock {\em International Journal for Numerical Methods in Fluids},
  63:313--340.

\bibitem[Legout et~al., 2012]{Legout12}
Legout, C., Darboux, F., N\'ed\'elec, Y., Hauet, A., Esteves, M., Renaux, B.,
  Denis, H., and Cordier, S. (2012).
\newblock High spatial resolution mapping of surface velocities and depths for
  shallow overland flow.
\newblock {\em Earth Surface Processes and Landforms}, 37(9):984--993.

\bibitem[Liang and Marche, 2009]{Liang09b}
Liang, Q. and Marche, F. (2009).
\newblock Numerical resolution of well-balanced shallow water equations with
  complex source terms.
\newblock {\em Advances in Water Resources}, 32(6):873--884.

\bibitem[Mac{D}onald, 1996]{MacDonald96}
Mac{D}onald, I. (1996).
\newblock {\em Analysis and computation of steady open channel flow}.
\newblock PhD thesis, University of Reading --- Department of Mathematics.

\bibitem[Mac{D}onald et~al., 1997]{MacDonald97}
Mac{D}onald, I., Baines, M.~J., Nichols, N.~K., and Samuels, P.~G. (1997).
\newblock Analytic benchmark solutions for open-channel flows.
\newblock {\em Journal of Hydraulic Engineering}, 123(11):1041--1045.

\bibitem[Mangeney et~al., 2007]{Mangeney07}
Mangeney, A., Bouchut, F., Thomas, N., Vilotte, J.-P., and Bristeau, M.-O.
  (2007).
\newblock Numerical modeling of self-channeling granular flows and of their
  levee-channel deposits.
\newblock {\em Journal of Geophysical Research}, 112(F2):F02017.

\bibitem[Mein and Larson, 1973]{Mein73}
Mein, R.~G. and Larson, C.~L. (1973).
\newblock Modeling infiltration during a steady rain.
\newblock {\em Water Resources Research}, 9(2):384--394.

\bibitem[Moussa and Bocquillon, 2000]{Moussa00}
Moussa, R. and Bocquillon, C. (2000).
\newblock Approximation zones of the {S}aint-{V}enant equations for flood
  routing with overbank flow.
\newblock {\em Hydrology and Earth System Sciences}, 4(2):251--260.

\bibitem[M\"ugler et~al., 2011]{Mugler11}
M\"ugler, C., Planchon, O., Patin, J., Weill, S., Silvera, N., Richard, P., and
  Mouche, E. (2011).
\newblock Comparison of roughness models to simulate overland flow and tracer
  transport experiments under simulated rainfall at plot scale.
\newblock {\em Journal of Hydrology}, 402(1--2):25--40.

\bibitem[Novak et~al., 2010]{Novak10}
Novak, P., Guinot, V., Jeffrey, A., and Reeve, D.~E. (2010).
\newblock {\em Hydraulic modelling --- An Introduction. Principles, Methods and
  Applications}.
\newblock Spon Press.

\bibitem[Paquier, 1995]{Paquier95}
Paquier, A. (1995).
\newblock {\em Modeling and simulation of the propagation of the dam-break
  wave}.
\newblock PhD thesis, University of Saint-Etienne, France.
\newblock In French.

\bibitem[Peng, 2011]{Peng11}
Peng, R.~D. (2011).
\newblock Reproducible research in computational science.
\newblock {\em Science}, 334(6060):1226--1227.

\bibitem[Popinet, 2011]{popinet11}
Popinet, S. (2011).
\newblock Quadtree-adaptative tsunami modelling.
\newblock {\em Ocean Dynamics}, 61(9):1261--1285.

\bibitem[Ritter, 1892]{Ritter92}
Ritter, A. (1892).
\newblock Die {F}ortpflanzung der {W}asserwellen.
\newblock {\em {Z}eitschrift des {V}ereines {D}euscher {I}ngenieure},
  36(33):947--954.

\bibitem[Savage and Hutter, 1991]{Savage91}
Savage, S. and Hutter, K. (1991).
\newblock The dynamics of avalanches of granular materials from initiation to
  runout. part i: Analysis.
\newblock {\em Acta Mechanica}, 86(1--4):201--223.

\bibitem[Stodden et~al., 2013]{Stodden13}
Stodden, V., Bailey, D.~H., Borwein, J., LeVeque, R.~J., Rider, W., and Stein,
  W. (2013).
\newblock Setting the default to reproducible: Reproducibility in computational
  and experimental mathematics.
\newblock Technical report, Institute for Computational and Experimental
  Research in Mathematics.

\bibitem[Tanguy and Chocat, 2013]{Tanguy13}
Tanguy, J.-M. and Chocat, B. (2013).
\newblock {\em CANOE}, pages 209--218.
\newblock John Wiley \& Sons, Inc.

\bibitem[Tatard et~al., 2008]{Tatard08}
Tatard, L., Planchon, O., Wainwright, J., Nord, G., Favis-Mortlock, D.,
  Silvera, N., Ribolzi, O., Esteves, M., and Huang, C. (2008).
\newblock Measurement and modelling of high-resolution flow-velocity data under
  simulated rainfall on a low-slope sandy soil.
\newblock {\em Journal of Hydrology}, 348(1-2):1--12.

\bibitem[Thacker, 1981]{Thacker81}
Thacker, W.~C. (1981).
\newblock Some exact solutions to the nonlinear shallow-water wave equations.
\newblock {\em Journal of Fluid Mechanics}, 107:499--508.

\bibitem[Toro et~al., 1994]{Toro94}
Toro, E., Spruce, M., and Speares, W. (1994).
\newblock Restoration of the contact surface in the {H}{L}{L}-{R}iemann solver.
\newblock {\em Shock Waves}, 4:25--34.

\bibitem[van Heesch, 2013]{doxygen}
van Heesch, D. (2013).
\newblock {\em {d}oxygen. Manual for version 1.8.5}.

\bibitem[Zhang and Cundy, 1989]{Zhang89}
Zhang, W. and Cundy, T.~W. (1989).
\newblock Modeling of two-dimensional overland flow.
\newblock {\em Water Resources Research}, 25(9):2019--2035.

\end{thebibliography}

\end{document}